\documentclass{SCAE}
\numberwithin{equation}{section}
\begin{document}

\Vol{} %
\No{} %
\BeginPage{1} %
\EndPage{XX} %
\AuthorMark{Fan X Q {\it et al.}}
\ReceivedDay{September 23, 2014}
\AcceptedDay{January  12, 2015}
\DOI{10.1007/s11425-000-0000-0} 

\title{Sharp large deviation results for sums of independent random
variables}{}


\author[1,2]{FAN Xiequan}{}
\author[3]{GRAMA Ion}{}
\author[3,4]{LIU Quansheng}{Corresponding author}

\address[{\rm1}]{Regularity Team, INRIA Saclay,  Palaiseau, 91120, France;   }%
\address[{\rm2}]{Department of Mathematics, School of Science, Tianjin University, Tianjin, 300072,   China;}
\address[{\rm3}]{ LMBA, UMR 6205, Univ. Bretagne-Sud,  Vannes, 56000, France;}
\address[{\rm4}]{College of Mathematics and Computing Science, Changsha University of Science and Technology, Changsha, 410004,  China}
\Emails{fanxiequan@hotmail.com,
ion.grama@univ-ubs.fr, quansheng.liu@univ-ubs.fr}\maketitle


 {\begin{center}
\parbox{14.5cm}{\begin{abstract}
We show sharp bounds for probabilities of large deviations for sums of independent random variables satisfying Bernstein's condition.  One such bound is very close to the tail of the standard Gaussian law in certain case; other bounds improve the inequalities of Bennett and Hoeffding  by adding missing factors in the spirit of Talagrand (1995). We also complete Talagrand's inequality by giving a lower bound of the same form, leading to an equality.  As a consequence, we obtain large deviation expansions similar to those of Cram\'{e}r (1938), Bahadur-Rao (1960) and Sakhanenko (1991). We also show that our bound can be used to improve a recent inequality of Pinelis (2014).\vspace{-3mm}
\end{abstract}}\end{center}}



\renewcommand{\baselinestretch}{1.2}
\begin{center} \renewcommand{\arraystretch}{1.5}
{\begin{tabular}{lp{0.8\textwidth}} \hline \scriptsize
{\bf Citation:}\!\!\!\!&\scriptsize Fan X Q, Grama I, Liu Q S.  Sharp large deviation results for sums of independent random
variables. Sci China Math,  2015, 58, doi: 10.1007/s11425-000-0000-0 \vspace{1mm}
\\
\hline
\end{tabular}}\end{center}

\baselineskip 11pt\parindent=10.8pt  \wuhao

\section{Introduction}
Let $\xi_{1},...,\xi _{n}$ be a finite sequence of independent centered random variables (r.v.). Denote by
\begin{equation}
S_{n}=\sum_{i=1}^{n}\xi _{i}\ \ \  \mbox{and}\ \ \
 \sigma ^{2}=\sum_{i=1}^{n}\mathbf{E}[\xi_i^2]. \label{SSI}
\end{equation}
Starting from the seminal work of Cram\'{e}r \cite{Cramer38}
and Bernstein \cite{B46}, the estimation of the tail
probabilities $\mathbf{P}\left( S_{n}>x\right) ,$ for large $x>0,$ has attracted much
attention. Various precise inequalities and asymptotic results have been established by Hoeffding \cite{Ho63},
Nagaev \cite{N79}, Saulis and  Statulevicius \cite{SS91},   Chaganty and Sethuraman \cite{CS93} and Petrov \cite{Petrov95}  under different backgrounds.

Assume that $(\xi_i)_{i=1,...,n}$ satisfies Bernstein's condition
\begin{equation}
|   \mathbf{E}[   \xi _{i}^{k}]|\leq \frac{1}{2}k!\varepsilon ^{k-2}\mathbf{E}[\xi _{i}^{2}],\ \ \ \mbox{for
}\ k\geq 3\ \ \mbox{and}\ \ i=1,...,n,  \label{Bernstein cond}
\end{equation}
for some constant $\varepsilon>0$.
By employing the exponential Markov
inequality and an upper bound for the moment generating function $\mathbf{E}[e^{\lambda \xi _{i}}]$,
Bernstein \cite{B46} (see also Bennett \cite{Be62})  has  obtained the following
inequalities: for all $x\geq0$,
\begin{eqnarray}
\mathbf{P}(S_{n}> x\sigma )&\leq& \inf_{\lambda\geq 0}\mathbf{E} [e^{\lambda(S_{n}- x\sigma)}]  \label{fg31} \\
&\leq& B\left( x, \frac \varepsilon \sigma  \right)
:=\exp \left\{ - \frac{\widehat{x}^2 }{2}  \right\}  \label{Beiq} \\
&\leq& \exp \left\{ - \frac{x^2 }{2 (1+ x \varepsilon/\sigma )}  \right\},  \label{Bersteidnsiq}
\end{eqnarray}
where
\begin{eqnarray}\label{Beiqt}
\widehat{x} =  \frac{2 x }{1+ \sqrt{1+2x\varepsilon/\sigma}};
\end{eqnarray}
see also  van de Geer and Lederer \cite{VL11} with a new method based  on Bernstein-Orlicz norm and   Rio \cite{R12}.
Some extensions of the inequalities of  Bernstein and Bennett can be found in van de Geer \cite{VL95}  and de la Pe\~{n}a \cite{De99} for martingales; see also  Rio \cite{R02,R12q} and Bousquet \cite{B02}  for the empirical processes with r.v. bounded from above.

Since $\lim_{\varepsilon /\sigma \rightarrow 0}\mathbf{P}(S_{n}> x\sigma )=1-\Phi(x)$ and $\lim_{\varepsilon /\sigma \rightarrow 0} B\left( x, \frac \varepsilon \sigma  \right)= e^{-x^2/2}$, where
\[
\Phi(x)=\frac{1}{\sqrt{2 \pi}}\int_{-\infty}^{x} e^{-\frac{t^2}{2}}dt
\]
is the standard normal distribution function, the central limit theorem (CLT) suggests that Bennett's inequality (\ref{Beiq}) can be substantially refined by adding the factor
\begin{eqnarray*}
M(x)=\Big(1-\Phi(x) \Big)\exp\left\{ \frac{x^2}{2}\right\},\ \label{mill}
\end{eqnarray*}
 where $\sqrt{2\pi}M(x)$ is known as Mill's ratio.

To recover the factor $M(x)$ of order $1/x$ as $x\rightarrow \infty$, a lot of effort has been made.
Certain factors of order $1/x$ have been recovered by using the following inequality: for some $\alpha>1,$
\begin{eqnarray*}
\mathbf{P}(S_{n}\geq x\sigma ) & \leq & \inf_{t< x\sigma} \mathbf{E}\bigg[ \frac{((S_{n}-t)^+)^\alpha}{((x\sigma-t)^+)^\alpha } \bigg],
\end{eqnarray*}
where $x^+= \max\{x, 0\}$; see  Eaton \cite{E74}, Bentkus \cite{Be03},  Pinelis \cite{P12} and Bentkus et al.\ \cite{BZ06}.
Some bounds on tail probabilities of type
\begin{eqnarray}\label{sfffs}
\mathbf{P}(S_{n}\geq x\sigma )\  \leq\  C \Big(1-\Phi(x) \Big) ,
\end{eqnarray}
where $C>1$ is an absolute constant, are obtained for sums of weighted Rademacher random variables; see Bentkus \cite{Be03}.
In particular, Bentkus and   Dzindzalieta \cite{BD14}  proved that $$C = \frac {1}{4 (1-\Phi(\sqrt{2}))}   \approx 3.178$$
is sharp in (\ref{sfffs}).

When the summands $\xi_{i}$ are bounded from above, results of such type have been obtained by Talagrand \cite{Ta95},  Bentkus \cite{Be04} and Pinelis \cite{P14}.
Using  the conjugate measure technique, Talagrand (cf.\ Theorems 1.1 and 3.3 of \cite{Ta95}) proved that if the random variables satisfy $\xi_{i}\leq 1$ and   $|\xi_i|\leq b$ for a constant $b>0$ and all $i=1,...,n,$ then there exists an universal constant $K$ such that, for all $0\leq x \leq  \frac{\sigma}{K b}$,
\begin{eqnarray}
\mathbf{P}(S_{n} > x\sigma) &\leq& \inf_{\lambda\geq 0}\mathbf{E}[e^{\lambda(S_n-x\sigma)} ]\left( M(x)+ K  \frac{b}{\sigma} \right)  \label{fjksa} \\
&\leq& H_n(x,\sigma) \left( M(x)+ K \frac{b}{\sigma} \right), \label{Ta}
\end{eqnarray}
where $$H_n(x,\sigma)=\left\{\left(\frac{\sigma}{x+\sigma}\right)^{x\sigma+\sigma^2}\left(\frac{n}{n-x\sigma}\right)^{n-x\sigma} \right\}^{\frac{n}{n+\sigma^2}}.$$
Since $M(x)= O\left( \frac{1}{x} \right), x \rightarrow \infty$,  equality (\ref{Ta}) improves on Hoeffding's bound $H_n(x,\sigma)$ (cf.\ (2.8) of \cite{Ho63}) by adding a missing factor $$F_1 (x, b/\sigma)=M(x)+ K b/ \sigma $$  of order $\frac{1}{x}$ for all
$0\leq x\leq \frac{\sigma}{K b} $.
 Other  improvements on  Hoeffding's bound  can be found in Bentkus \cite{Be04} and Pinelis \cite{P12}.
Bentkus's inequality \cite{Be04} is much better than (\ref{Ta}) in the sense that it recovers a factor of order $\frac{1}{x}$ for all $x\geq 0$ instead of the range $0\leq x\leq \frac{\sigma}{K b}$, and do not assume that $\xi_i$'s have moments of order larger than $2$; see also
Pinelis \cite{P14} for a similar improvement on Bennett-Hoeffding's bound.

The scope of this paper is to give several improvements on  Bernstein's inequalities (\ref{fg31}), (\ref{Bersteidnsiq})  and Bennett's inequality  (\ref{Beiq}) for sums of non-bounded random variables instead of sums of bounded (from above) random variables, which are considered in Talagrand \cite{Ta95}, Bentkus \cite{Be04} and Pinelis \cite{P12}.  Moreover, some tight lower bounds are also given, which were not considered by Talagrand \cite{Ta95}, Bentkus \cite{Be04} and Pinelis \cite{P12}. In particular, we improve Talagrand's inequality to an {\em equality}, which will imply simple large deviation expansions.
 We also show that our bound can be used to improve a recent upper bound on tail probabilities due to Pinelis \cite{P12}.

Our approach is based on the conjugate distribution technique due to Cram\'{e}r, which becomes a standard for obtaining sharp large deviation expansions. We refine the technique inspired by Talagrand \cite{Ta95} and  Grama and Haeusler \cite{GH00} (see also \cite{FGL13,F13}), and derive sharp bounds for the cumulant function to obtain  precise upper bounds on tail probabilities under Bernstein's condition.

As to the potential applications of our results in statistics, we refer to Fu, Li and Zhao \cite{F93} for large sample estimation  and Joutard \cite{J06,J12} for nonparametric estimation. In these papers, many interesting Bahadur-Rao type large deviation expansions have been established. Our result leads to simple large deviation expansions which  are similar (but simpler) to those  of Cram\'{e}r (1938), Bahadur-Rao (1960) and Sakhanenko (1991).     For other important applications, we refer to Shao \cite{S99} and Jing, Shao and Wang \cite{JSW03}, where the authors have established the Cram\'{e}r type self-normalized large deviations for normalized  $x =o(n^{1/6});$ see also Jing, Liang and Zhou \cite{JLZ12}. From the proofs of theorems in \cite{S99,JSW03,JLZ12}, we find that the self-normalized large deviations are closely related  to the large deviations for sums of bounded from above random variables (cf.\ \cite{F12}). Our results may help extend the Cram\'{e}r type self-normalized  large deviations to  a larger range.

The paper is organized as follows. In Section \ref{sec2}, we present our main
results. In Section \ref{secco},  some important comparisons are given.
In Section \ref{sec2.5}, we state some auxiliary results to be used in the proofs of theorems.
Sections \ref{sec3.5} - \ref{secend} are devoted to the proofs of main results.

\section{Main results}\label{sec2}
All over the paper $\xi_{1}, ..., \xi_{n}$ is  a finite sequence of independent real random variables with $\mathbf{E}[\xi_{i}]=0$ and satisfying Bernstein's condition (\ref{Bernstein cond}), $S_n$ and $\sigma^2$
are defined by (\ref{SSI}). We use the notations $a\wedge b=\min\{a, b\}$, $a\vee b=\max\{a, b\}$ and $a^{+}=a\vee 0$.
Throughout this paper, $C$ stands for an absolute constant with possibly different values in different places.

Our first result is the following large deviation inequality valid  for all $x\geq0$.
\begin{theorem}
\label{th2}For any $\delta\in(0,1]$ and $x\geq 0$,
\begin{eqnarray}  \label{f10}
\ \ \ \ \ \ \mathbf{P}(S_n>x\sigma ) &\leq& \Big(1-\Phi\left(  \widetilde{x}\right) \Big)  \left[1+C_\delta \left(1+ \widetilde{x}\right) \frac{\varepsilon}{\sigma} \right],
\end{eqnarray}
where
\begin{eqnarray}  \label{f1dgfr}
\widetilde{x}=\frac{2x}{1+\sqrt{1+2(1+\delta)x\varepsilon/\sigma}}
\end{eqnarray}
and $C_\delta$ is a constant only depending on $\delta$.  In particular, if $0\leq x = o(\sigma/\varepsilon), \  \varepsilon/\sigma \rightarrow 0,$ then
\begin{eqnarray*}
\ \ \ \ \ \ \mathbf{P}(S_n>x\sigma ) &\leq& \Big(1-\Phi\left(  \widetilde{x}\right) \Big)  \Big[1+o(1)\Big].
\end{eqnarray*}
\end{theorem}

The interesting feature of the bound (\ref{f10}) is that it decays exponentially to $0$ and also recovers closely the shape of the standard normal tail $1-\Phi(x)$ when  $r= \frac \varepsilon \sigma$ becomes small, which is not the case of Bennett's bound $B(x, \frac \varepsilon \sigma)$ and Berry-Essen's bound
\[
\mathbf{P}(S_n>x\sigma )   \leq  1-\Phi\left( x \right) +  C \frac{\varepsilon}{\sigma}.
\]
Our result can  be compared with Cram\'{e}r's large deviation result in the i.i.d.\ case (cf.\ (\ref{cramesrse})). With respect to Cram\'{e}r's result, the advantage of (\ref{f10}) is that it is valid for all $x\geq0$.

Notice that Theorem \ref{th2} improves Bennett's  bound only for moderate $x$.
A further significant improvement of Bennett's inequality (\ref{Beiq}) for all $x\geq 0$ is given by the following theorem: we replace Bennett's bound $B\left( x, \frac \varepsilon \sigma \right) $ by
the following smaller one:
\begin{eqnarray}
  B_{n}\left(x, \frac{\varepsilon}{\sigma} \right) &=& B\left(x, \frac{\varepsilon}{\sigma} \right) \exp\left\{-n \psi\left(\frac{\widehat{x}^2}{2n\sqrt{1+2x \varepsilon/\sigma}} \right)\right\},
\end{eqnarray}
where $\psi(t)=t-\log(1+t)$ is a nonnegative convex function in $t\geq 0$.

\begin{theorem}
\label{th3} For all $x\geq 0$,
\begin{eqnarray}
 \mathbf{P}(S_{n} > x\sigma ) &\leq& B_{n}\left(x, \frac{\varepsilon}{\sigma} \right)\, F_{2}\left(x, \frac{\varepsilon}{\sigma}\right)  \label{f17}\\
 &\leq& B_{n}\left(x, \frac{\varepsilon}{\sigma} \right),
\end{eqnarray}
where
\begin{eqnarray}
F_{2}\left(x, \frac{\varepsilon}{\sigma}\right) =   \left( M(x) +  27.99 R\left( x \varepsilon/\sigma  \right) \frac{\varepsilon}{\sigma} \right)\wedge 1
\end{eqnarray}
and
\begin{eqnarray}\label{funR}
R(t) = \left\{ \begin{array}{ll}
\frac{(1-t+6t^2)^3 }{(1-3t)^{3/2}(1-t)^7 },\ \  & \textrm{if $0\leq t< \frac13$},\\
\infty, & \textrm{if $t\geq \frac13$,}
\end{array} \right.
\end{eqnarray}
is an increasing function. Moreover,  for all $0\leq x \leq \alpha \frac{\sigma}{\varepsilon}$ with $0\leq \alpha < \frac 1 3$, it holds $R( x \varepsilon/\sigma  )\leq R(\alpha)$.
 If $\alpha=0.1$, then $27.99R(\alpha) \leq  88.41$.
\end{theorem}

To highlight the improvement of Theorem \ref{th3} over Bennett's bound, we note that
$B_{n}(x, \frac \varepsilon \sigma )\leq B(x, \frac \varepsilon  \sigma )$ and, in the i.i.d.\ case (or, more generally when $\frac{\varepsilon}{\sigma}=\frac{c_0}{\sqrt{n}}$, for some constant $c_0>0$),
\begin{eqnarray}\label{fbnb}
 \ \ \ \ \ \ \ \ \  B_{n}\left( \sqrt{n}x, \frac{\varepsilon}{\sigma} \right)  &=&B\left( \sqrt{n}x, \frac{\varepsilon}{\sigma} \right)  \exp \left\{  -  c_x \, n \right\},
\end{eqnarray}
where $c_x>0$, $x>0$, does not depend on $n$. Thus Bennett's bound is strengthened by
adding a factor $\exp \left\{  -  c_x \, n \right\}, n \rightarrow \infty,$ which is similar to Hoeffding's improvement on Bennett's bound for sums of bounded random variables \cite{Ho63}.
\begin{figure}\center
\includegraphics[width=0.4\textwidth]{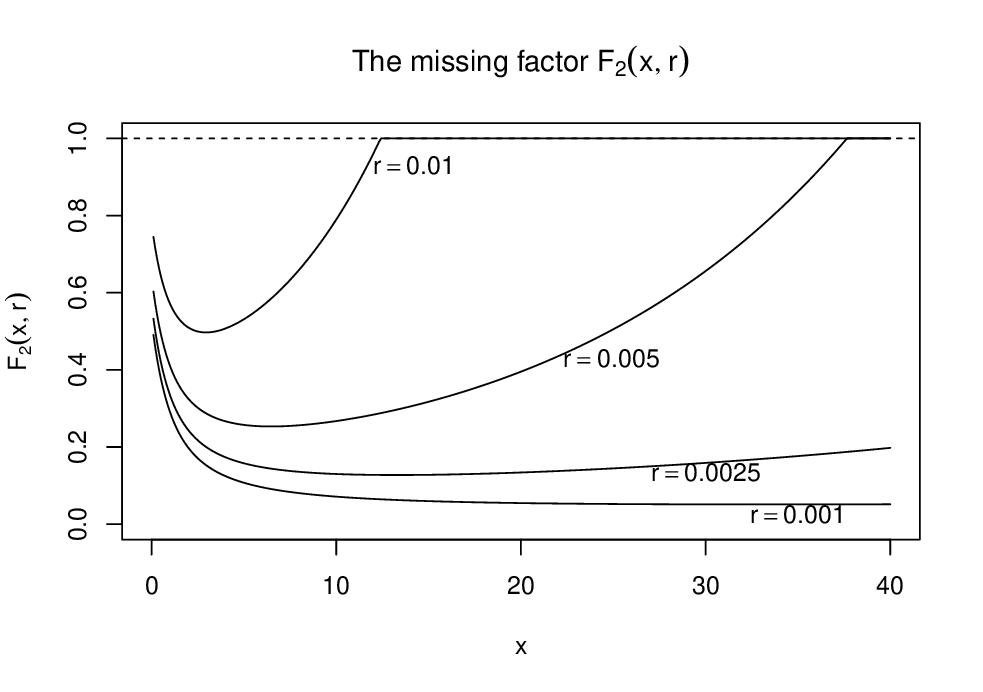}
\caption[]{ The missing factor $F_2(x, r)$ is displayed as a function of $x$ for various values of $r=\frac{\varepsilon}{\sigma}$.}
\label{Th4a}
\end{figure}
The second improvement in the right-hand side of (\ref{f17}) comes from the missing factor $F_{2}(x, \frac \varepsilon \sigma )$, which is of order $M(x)[1+o(1)]$ for moderate values of $x$ satisfying
$0\leq x= o( \frac{\sigma}{ \varepsilon}),  \ \frac{ \varepsilon}{\sigma}\rightarrow 0$. This improvement is similar  to Talagrand's refinement on Hoeffding's upper bound $H_n(x,\sigma)$ by the factor $F_1(x,b/\sigma)$; see
 (\ref{Ta}). The numerical values of the missing factor $F_2(x, \frac \varepsilon \sigma)$
are displayed in Figure \ref{Th4a}.
\begin{figure}\center
\includegraphics[width=0.4\textwidth]{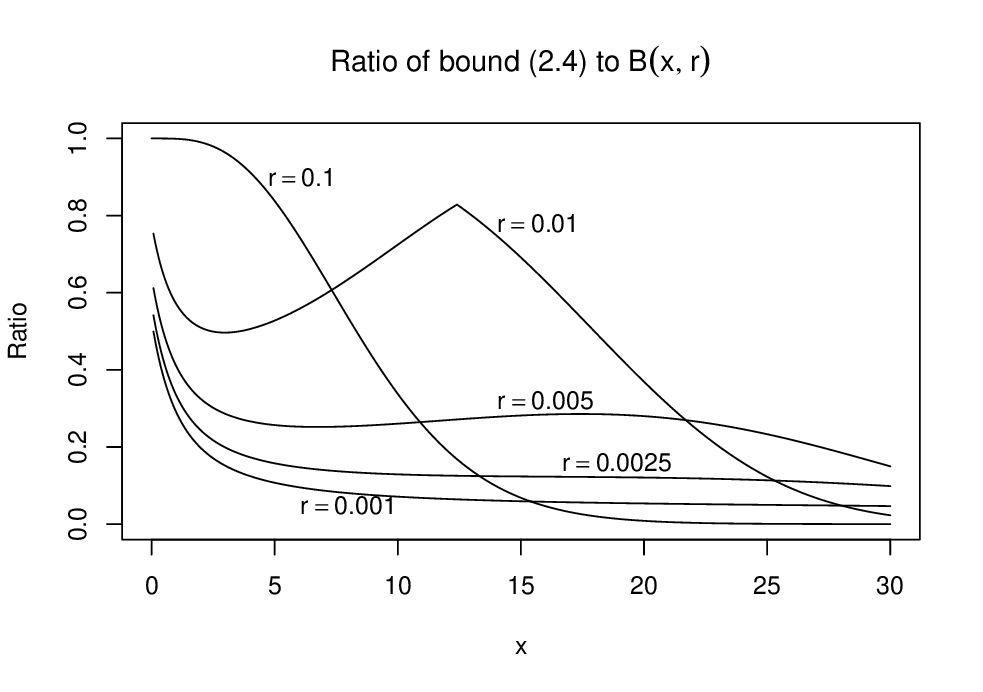}
\caption[]{ Ratio of $B_{n}\left(x, r \right) F_{2}\left(x, r\right)$ to  $B(x, r)$ as a function of $x$ for various values of $r=\frac{\varepsilon}{\sigma}=\frac{1}{\sqrt{n}}$.}
\label{Th3a}
\end{figure}

Our numerical results confirm that the bound $B_{n}(x, \frac \varepsilon \sigma) F_{2}(x,\frac \varepsilon \sigma)$ in (\ref{f17}) is
better than Bennett's bound $B(x,\frac \varepsilon \sigma)$ for all $x\geq 0$. For the convenience of the reader, we display
the ratios of $B_n(x, r)F_2(x, r)$ to $B(x, r)$ in Figure~\ref{Th3a} for various $r= \frac{1}{\sqrt{n}}$.

The following corollary improves inequality (\ref{f10}) of Theorem \ref{th2} in the range $0\leq x \leq \alpha \frac \sigma \varepsilon$ with $0\leq\alpha< \frac 13$. It corresponds to taking $\delta=0$ in the definition (\ref{f1dgfr}) of $\widetilde{x}$.
\begin{corollary}\label{co4}For all $0\leq x \leq \alpha  \frac \sigma \varepsilon$ with $0\leq \alpha < \frac 13$,
\begin{eqnarray}
\ \ \ \ \ \ \mathbf{P}(S_{n} > x\sigma )&\leq& \Big(1-\Phi \left(\widehat{x} \right) \Big) \left[  1 +  70.17 R(\alpha)  \left(1+\widehat{x} \right) \,\frac{\varepsilon  }{ \sigma }\right] , \label{fmb}
\end{eqnarray}
where $\widehat{x}$ is defined in (\ref{Beiqt}) and $R(t)$ by (\ref{funR}). In particular, for all $0\leq x = o(\frac \sigma \varepsilon), \ \frac \varepsilon\sigma \rightarrow 0$,
\begin{eqnarray}
\ \ \ \ \ \ \mathbf{P}(S_{n} > x\sigma )&\leq& \Big(1-\Phi \left(\widehat{x} \right) \Big) \Big[  1 +  o(1)\Big]\label{sxkmb}\\
&=& B\left(x, \frac \varepsilon \sigma \right) M(\widehat{x})\Big[ 1+ o(1) \Big] . \nonumber
\end{eqnarray}
\end{corollary}

The advantage of Corollary \ref{co4} is that in the normal distribution function
$\Phi(x)$ we have the expression $\widehat{x}$ instead of the smaller term $\widetilde{x}$ figuring in Theorem \ref{th2},
which represents a significant improvement.

Notice that inequality (\ref{sxkmb}) improves Bennett's bound $B\left(x, \frac \varepsilon \sigma \right)$  by the missing factor $M(\widehat{x})[ 1+ o(1) ]$ for all $0\leq x =o\left(\frac \sigma \varepsilon  \right)$.

For the lower bound of tail probabilities $\mathbf{P}(S_n>x\sigma )$, we have the
following result, which is a complement of Corollary \ref{co4}.

\begin{theorem}
\label{th4} For all $0\leq x \leq  \alpha  \frac \sigma  \varepsilon $ with $0\leq \alpha\leq \frac 1 {9.6}$,
\begin{eqnarray}\label{f19}
\mathbf{P}(S_n>x\sigma ) &\geq& \Big( 1-\Phi\left( \check{x} \right) \Big)  \left[1 - c_\alpha \left(1+\check{x} \right)\frac{\varepsilon }{\sigma}\right],\nonumber
\end{eqnarray}
where $\check{x}=\frac{ \lambda  \sigma}{(1- \lambda \varepsilon)^3}$ with $\lambda =\frac{2x/\sigma }{1 +\sqrt{1-9.6 x \varepsilon /\sigma }}$, and
$c_\alpha=67.38R\left( \frac{2\alpha}{1+\sqrt{1-9.6 \alpha}} \right)$ is a bounded function. Moreover,  for all $0\leq x= o(\frac{\sigma}{\varepsilon}),$  \ $\frac{\varepsilon}{\sigma}\rightarrow 0,$
\[
\mathbf{P}(S_n>x\sigma )  \geq  \Big( 1-\Phi\left( \check{x} \right) \Big)  \Big[1 - o(1)\Big].
\]
\end{theorem}

Combining Corollary \ref{co4} and Theorem \ref{th4}, we obtain,
for all $0\leq x \leq  0.1 \frac \sigma  \varepsilon $,
\begin{eqnarray}\label{cct}
\ \ \ \ \ \ \mathbf{P}(S_n>x\sigma )=\Bigg(1-\Phi\Big(x(1+\theta_1 c_1 x \frac{\varepsilon}{\sigma})\Big)\Bigg)
\Bigg[1+ \theta_2 c_{2}(1+x)\frac{\varepsilon}{\sigma} \Bigg],
\end{eqnarray}
where $c_1, c_2>0$ are absolute constants and  $|\theta_1|, |\theta_2|\leq1$. This result can be found in Sakhanenko \cite{S91} but in a more narrow zone.

Some earlier lower bounds on tail probabilities, based on Cram\'{e}r large deviations, can be found in Arkhangelskii \cite{A89},  Nagaev \cite{N02} and Rozovky \cite{L03}. In particular, Nagaev established the following lower bound
\begin{eqnarray}\label{cfbd}
\mathbf{P}(S_n>x\sigma ) &\geq& \Big( 1-\Phi(x) \Big) e^{-c_1x^3\frac{\varepsilon}{\sigma} } \left( 1- c_2 (1+x)\frac{\varepsilon}{\sigma} \right)
\end{eqnarray}
for some explicit constants $c_1, c_2$ and all $0 \leq x \leq \frac{1}{25} \frac{\sigma}{\varepsilon}$. For more general results, we refer to Theorem 3.1 of Saulis and  Statulevicius \cite{SS91}.

In the following theorem, we obtain a one term sharp large deviation expansion similar to Cram\'{e}r \cite{Cramer38},  Bahadur-Rao \cite{BR60}, Saulis and  Statulevicius \cite{SS91} and  Sakhanenko \cite{S91}.
\begin{theorem}
\label{thend} For all $0\leq x < \frac 1{12} \frac \sigma \varepsilon$,
\begin{eqnarray} \label{fed}
 \mathbf{P}(S_{n} > x\sigma ) = \inf_{\lambda\geq 0}\mathbf{E}[ e^{\lambda(S_n-x\sigma)} ] \, F_{3}\left(x, \frac{\varepsilon}{\sigma}\right) ,
\end{eqnarray}
where
\begin{eqnarray}
F_{3}\left(x, \frac{\varepsilon}{\sigma}\right) =   M(x) +  27.99\, \theta R\left(4 x\varepsilon/\sigma  \right) \frac{\varepsilon}{\sigma},\label{f3x}
\end{eqnarray}
$|\theta|\leq1$ and $R(t)$ is defined by (\ref{funR}). Moreover, $\inf_{\lambda\geq 0} \mathbf{E}[e^{\lambda(S_n-x\sigma)}] \leq B(x, \frac \varepsilon \sigma )$.  In particular, in the i.i.d.\ case, we have the following  non-uniform Berry-Esseen type bound: for all $0\leq x =o(\sqrt{n})$,
\begin{eqnarray}\label{jnskm}
\Big| \mathbf{P}(S_{n} > x\sigma )- M(x) \inf_{\lambda\geq 0} \mathbf{E}[e^{\lambda(S_n-x\sigma)}]\Big| \leq    \frac{C}{ \sqrt{n}}B(x, \frac \varepsilon \sigma )  .
\end{eqnarray}
\end{theorem}

Theorem \ref{thend} holds also for $\xi_i$'s bounded from above.  In this case the term $27.99\, \theta R\left(4 x\varepsilon/\sigma  \right)$ can be significantly refined; see \cite{F12}. In particular, if $|\xi_i|\leq \varepsilon,$ then  $27.99\, \theta R\left(4 x\varepsilon/\sigma  \right)$ can be improved to $3.08.$ However, under the stated condition of Theorem \ref{thend}, the term $27.99\, \theta R\left(4 x\varepsilon/\sigma  \right)$ cannot be   improved significantly.

When Bernstein's condition  fails, we refer to Theorem 3.1 of Saulis and  Statulevicius \cite{SS91},
where  explicit and asymptotic expansions have been established  via  the Cram\'{e}r series (cf.\ Petrov \cite{Petrov75}    for details).
 When the Bernstein  condition holds, their result reduces to the
result of Cram\'{e}r \cite{Cramer38}. However, they gave an explicit information on the term corresponding to our term  $27.99\, \theta R\left(4 x\varepsilon/\sigma  \right).$

Equality (\ref{fed}) shows that $\inf_{\lambda\geq 0}\mathbf{E}[e^{\lambda(S_n-x\sigma)}]$ is the best possible exponentially decreasing rate on tail probabilities.
It reveals the missing factor $F_3$ in Bernstein's bound  (\ref{fg31}) (and thus in many other classical bounds such as Hoeffding, Bennett and Bernstein).
Since $\theta \geq -1$, equality (\ref{fed}) completes Talagrand's upper bound (\ref{fjksa}) by giving a sharp lower bound.
If $\xi_i$ are bounded from above $\xi_i \leq 1$, it holds that $\inf_{\lambda\geq0}\mathbf{E}[e^{\lambda(S_n-x\sigma)} ] \leq H_{n}(x, \sigma)$ (cf.\ \cite{Ho63}). Therefore (\ref{fed}) implies Talagrand's inequality (\ref{Ta}).

A   precise large deviation expansion, as sharp as (\ref{fed}), can be found in Sakhanenko \cite{S91} (see also Gy\"{o}rfi,  Harrem\"{o}es and Tusn\'{a}dy \cite{GHT12}).
In his paper, Sakhanenko proved an equality similar to (\ref{fed}) in a more narrow range $0\leq x \leq \frac 1{200} \frac{\sigma}{\varepsilon},$
\begin{eqnarray}
 \Big| \mathbf{P}(S_{n} > x\sigma ) - \Big( 1- \Phi(t_x) \Big)   \Big| \ \leq \  C \frac{\varepsilon}{\sigma} e^{-t_x^2 /2}  ,
\end{eqnarray}
 where
$$t_x= \sqrt{-2 \ln \Big( \inf_{ \lambda \geq 0}\mathbf{E}[e^{\lambda(S_n-x\sigma)} ]\Big) } $$
is a value depending on the distribution of $S_n$ and satisfying $|t_x- x| =O(x^2 \frac{\varepsilon}{\sigma}), \ \frac{\varepsilon}{\sigma}  \rightarrow 0,$ for moderate $x$'s. It is worth noting that from Sakhanenko's result, we find that the inequalities (\ref{jnskm}) and (\ref{thols})   hold  also if $M(x)$ is replaced by $M(t_x)$.

\begin{figure}\center
\includegraphics[width=0.49\textwidth]{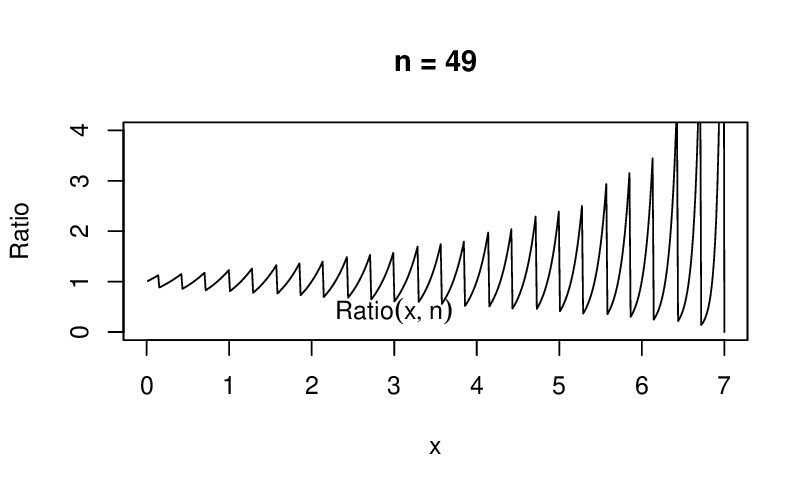}
\includegraphics[width=0.49\textwidth]{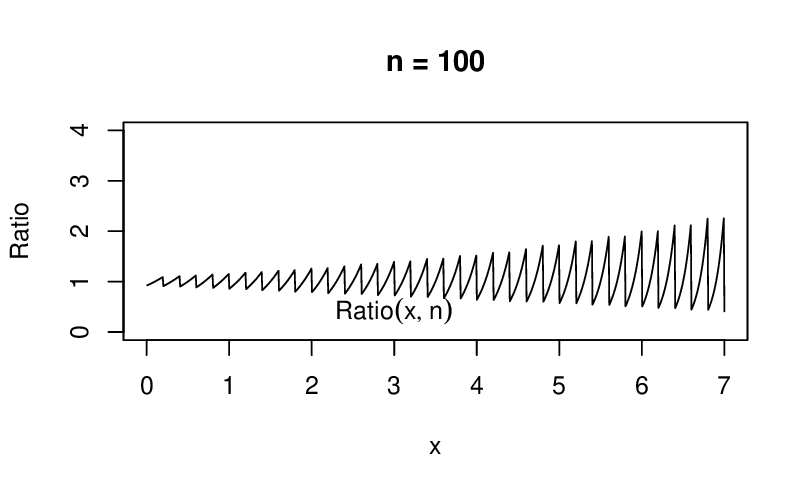}
\includegraphics[width=0.49\textwidth]{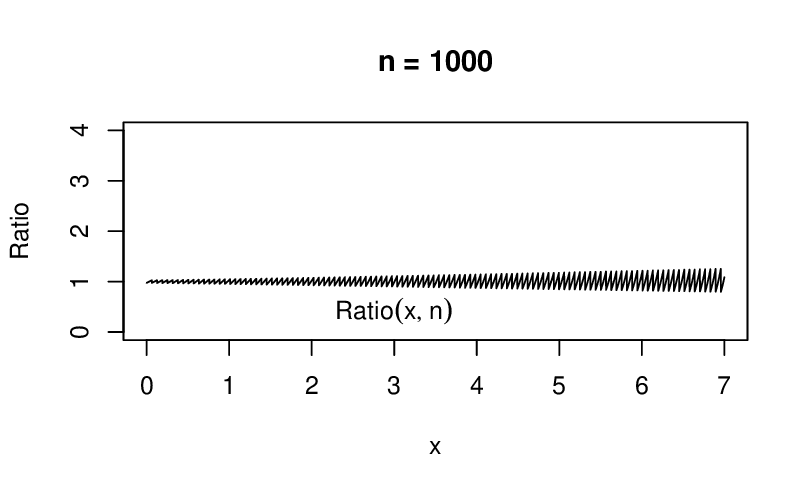}
\caption[]{The ratio  $\textrm{Ratio}(x, n)=\frac{\textbf{\textrm{P}}(S_{n} \geq x \sqrt{n})}{M(x)\inf_{\lambda \geq 0}\textbf{\textrm{E}}[e^{\lambda(S_n-x\sqrt{n})}]}$ is displayed as a function of $x$ for various $n$ for sums of Rademacher random variables.  }\label{Sim}
\end{figure}

Using the two sided bound
\begin{eqnarray}\label{mt}
 \frac{1}{\sqrt{2\pi}(1+t)} \leq M(t)\leq \frac{1}{\sqrt{\pi}(1+t)}, \ \ \ \ \ t\geq0,
\end{eqnarray}
and
\begin{eqnarray*}
 M(t) \sim \frac{1}{\sqrt{2\pi}(1+t)} , \ \  \ \ \ \ t \rightarrow \infty
\end{eqnarray*}
(see p.\ 17 in It\={o} and MacKean  \cite{IM96} or Talagrand \cite{Ta95}), equality (\ref{fed}) implies that the relative errors  between $\mathbf{P}(S_{n} > x\sigma )$ and
$M(x)\inf_{\lambda\geq 0} \mathbf{E}[e^{\lambda(S_n-x\sigma)}]$ converges to $0$ uniformly in the range $0\leq x =o\left( \frac{\sigma}{\varepsilon}\right)$ as $\frac{\varepsilon}{\sigma}\rightarrow 0$, i.e.
\begin{eqnarray} \label{thols}
\mathbf{P}(S_{n} > x\sigma )
 &=& M(x)\inf_{\lambda\geq 0} \mathbf{E}[e^{\lambda(S_n-x\sigma)}] \Big( 1+ o(1) \Big) .
\end{eqnarray}
Expansion (\ref{thols}) extends the following Cram\'{e}r large deviation expansion: for $0\leq x =o\left( \sqrt[3]{\frac{\sigma}{\varepsilon}}\right)$ as $\frac{\sigma}{\varepsilon}\rightarrow \infty$,
\begin{eqnarray} \label{tholv}
\ \ \ \ \ \ \mathbf{P}(S_n>x\sigma ) &=& \Big(1-\Phi\left( x\right) \Big) \Big[1+o(1)\Big].
\end{eqnarray}
To have an idea of the precision of expansion (\ref{thols}), we plot the ratio
$$\textrm{Ratio}(x, n)=\frac{\textbf{\textrm{P}}(S_{n} \geq x \sqrt{n})}{M(x)\inf_{\lambda \geq 0}\textbf{\textrm{E}}[e^{\lambda(S_n-x\sqrt{n})}]}$$
in Figure \ref{Sim} for the case of sums of Rademacher random variables $\mathbf{P}(\xi_i=-1)=\mathbf{P}(\xi_i=1)=\frac{1}{2}$. From these plots
we see that the error in (\ref{thols}) becomes smaller as $n$ increases.

\section{Some comparisons}\label{secco}
\subsection{Comparison with a recent inequality  of Pinelis }
In this subsection, we show that Theorem \ref{thend}  can be used to improve a recent upper bound on tail probabilities due to Pinelis \cite{P14}.
For simplicity of notations, we assume that $\xi_i \leq 1$ and only consider the i.i.d.\ case. For other cases, the argument is similar. Let us recall the notations of Pinelis.
 Denote by $\Gamma_{a^2}$ the normal random variable with  mean $0$ and variance $a^2 >0$, and $\Pi_\theta$
the Poisson random variable with parameter $\theta > 0$. Let also
\[
\widetilde{\Pi}_\theta \sim \Pi_\theta -\theta.
\]
Denote by
\begin{eqnarray}   \label{hjks}
\delta  = \frac{ \sum_{i=1}^{n} \mathbf{E}[(\xi_i^+)^3]}{\sigma^2 } .
\end{eqnarray}
Then it is obvious that  $\delta \in (0, 1).$
Pinelis (cf.\ Corollary 2.2 of \cite{P14}) proved that: for all $y \geq 0,$
\begin{eqnarray}\label{pines}
\mathbf{P}(S_n > y)\ \leq \ \frac{2e^3 }{9 }  \mathbf{P}^{LC}( \Gamma_{(1-\delta)\sigma^2} +  \widetilde{\Pi}_{\delta \sigma^2 } > y),
\end{eqnarray}
where, for any r.v. $\zeta$, the function $\mathbf{P}^{LC}(\zeta > y)$ denotes the least log-concave
majorant of the tail function $\mathbf{P}(\zeta > y).$ So that $\mathbf{P}^{LC}(\zeta > y)\geq \mathbf{P}(\zeta > y)$. By the remark of Pinelis, inequality (\ref{pines}) refines the
Bennet-Hoeffding inequality by adding a factor of order $\frac1x$ in certain range.
By Theorem \ref{thend} and some simple calculations, we find that, for all $0\leq y =o(n)$,
\begin{eqnarray}
&& \mathbf{P}^{LC} ( \Gamma_{(1-\delta)\sigma^2} +  \widetilde{\Pi}_{\delta \sigma^2 } > y) \nonumber \\
&\geq& \mathbf{P} ( \Gamma_{(1-\delta)\sigma^2} +  \widetilde{\Pi}_{\delta \sigma^2 } > y) \nonumber \\
 \ &\geq& \    \inf_{\lambda\geq 0} \mathbf{E}[e^{\lambda(\Gamma_{(1-\delta)\sigma^2} +  \widetilde{\Pi}_{\delta \sigma^2 } -y)} ]\Big(  M( y/ \sigma) -   \frac{C}{\sqrt{n}}    \Big)  \nonumber\\
\ &=& \    \inf_{\lambda\geq 0} \mathbf{E}[e^{-\lambda y + f(\lambda,\delta,\sigma)  }] \Big(  M( y/ \sigma) - \frac{C}{\sqrt{n}}     \Big),  \label{sdxpinesli}
\end{eqnarray}
where
\[
f(\lambda,\delta,\sigma)=\frac{\lambda^2}{2}(1-\delta)\sigma^2 +  (e^\lambda -1 -\lambda )\delta \sigma^2 \, .
\]
By the inequality
\[
e^x-1-x \leq \frac{x^2}{2}  + \sum_{k=3}^{\infty}(x^+)^k
\]
(cf.\ proof of Corollary 3 in Rio \cite{R12}) and the fact that $\log(1+t)$ is concave in $t\geq 0$,
it follows that, for any $\lambda>0,$
\begin{eqnarray}
 \log \mathbf{E}[e^{\lambda S_n  }]
\ &=&  \sum_{i=1}^n \log \mathbf{E}[e^{\lambda \xi_i  } ]  \leq  \sum_{i=1}^n \log \Big(1+ \frac{\lambda^2}{2} \mathbf{E}[\xi_i^2 ] + \sum_{k=3}^{\infty}  \frac{\lambda^k}{k!} \mathbf{E}[(\xi_i^+)^k] \Big)  \nonumber\\
\ &\leq& \sum_{i=1}^n \log \Big(1+ \frac{\lambda^2}{2} \mathbf{E}[\xi_i^2] + \sum_{k=3}^{\infty}  \frac{\lambda^k}{k!} \mathbf{E}[(\xi_i^+)^3] \Big)  \nonumber\\
\ &\leq& n \log \Bigg(1+ \frac1n \sum_{i=1}^n \Big(\frac{\lambda^2}{2} \mathbf{E}[\xi_i^2] + \sum_{k=3}^{\infty}  \frac{\lambda^k}{k!} \mathbf{E}[(\xi_i^+)^3] \Big) \Bigg)  \nonumber\\
\ &\leq& n \log \Bigg(1+  \frac1n\Big(\frac{\lambda^2}{2}\sigma^2 +  \sum_{k=3}^{\infty}  \frac{\lambda^k}{k!} \delta \sigma^2\Big) \Bigg)  \nonumber\\
\ &=&  n \log \bigg(1+  \frac1nf(\lambda,\delta,\sigma) \bigg) .\nonumber
\end{eqnarray}
By the last  line, Theorem \ref{thend} implies that,
for all $0\leq y =o(n)$,
\begin{eqnarray}
\mathbf{P}(S_n > y)\ &\leq&  \  \inf_{\lambda\geq 0} \mathbf{E}[e^{-\lambda y + n \log (1+  \frac1nf(\lambda,\delta,\sigma) )   }] \Big(  M( y/ \sigma) + \frac{C}{\sqrt{n}}     \Big) \nonumber\\
&\leq&  \Big( 1+o(1) \Big) \inf_{\lambda\geq 0} \mathbf{E}[e^{-\lambda y + n \log (1+  \frac1nf(\lambda,\delta,\sigma) )   }] \Big(  M( y/ \sigma) - \frac{C}{\sqrt{n}}     \Big) .\label{sfc}
\end{eqnarray}
Note that $n \log \big(1+  \frac1nf(\lambda,\delta,\sigma) \big) \leq f(\lambda,\delta,\sigma).$
By the  inequalities (\ref{pines}), (\ref{sdxpinesli}) and (\ref{sfc}),
we find that (\ref{sfc}) not only refines  Pinelis' constant $\frac{2e^3 }{9 }\ (\approx 4.463)$ to $1+o(1)$ for large $n$, but also gives an exponential bound sharper than that of Pinelis.

\subsection{Comparison with the expansions of Cram\'{e}r and Bahadur-Rao }
Notice that the expression $\inf_{ \lambda \geq 0}\mathbf{E}[e^{\lambda(S_n-x\sigma)}]$ can be rewritten in the form $\exp\{-n \Lambda_n^\ast(\frac{x\sigma}n)\}$, where $\Lambda_n^\ast(x)=\sup_{\lambda\geq 0}\{\lambda x - \frac1n \log \mathbf{E}[e^{\lambda S_n}]\}$ is the Fenchel-Legendre transform  of the normalized cumulant  function of $S_n$. In the i.i.d.\ case, the function $\Lambda^\ast(x)= \Lambda_n^\ast(x)$ is known as the good rate function in large deviation principle (LDP) theory (see  Deuschel and Stroock \cite{D89} or Dembo and Zeitouni \cite{D98}).

Now we clarify the relation among our large deviation expansion (\ref{fed}), Cram\'{e}r large deviations \cite{Cramer38} and the Bahadur-Rao theorem \cite{BR60} in the i.i.d.\ case. Without loss of generality, we take $\sigma_1^2=1$, where $\sigma_1^2$
is the variance of $\xi_1$. First, our bound (\ref{fed}) implies that: for all $0\leq x =o(\sqrt{n})$,
\begin{eqnarray}\label{dsvss}
\mathbf{P}\left( S_{n}  > x \sigma \right)=   e^{   - n \, \Lambda^\ast(x/\sqrt{n})  }M(x) \left[ 1 + O\left(  \frac{1+x}{ \sqrt{n}}\right) \right],\ \ \ n\rightarrow \infty.
\end{eqnarray}
Cram\'{e}r \cite{Cramer38} (see also Theorem 3.1 of Saulis and   Statulevicius \cite{SS91} for more general results) proved that, for all $0\leq x =o(\sqrt{n})$,
\begin{equation}\label{cramesrse}
\frac{\mathbf{P}(S_{n}> x \sigma)}{1-\Phi(x)} =\exp\left\{ \frac{x^3}{\sqrt{n}}  \lambda\left(\frac{x}{\sqrt{n}} \right)\right\} \left[ 1+ O\left(  \frac{1+x}{ \sqrt{n}}\right) \right],\ \ n\rightarrow \infty,
\end{equation}
where $\lambda(\cdot)$ is   the Cram\'{e}r series.
So the good rate function and the Cram\'{e}r series   have the relation $n \, \Lambda^\ast(\frac{x}{\sqrt{n}})= \frac{x^2} 2 -\frac{x^3}{\sqrt{n}}  \lambda\left(\frac{x}{\sqrt{n}} \right).$
Second, consider the large deviation probabilities $\mathbf{P}\left( \frac{S_{n}}{n} > y \right)$.  Since $\frac{S_{n}}{n}\rightarrow 0, a.s.,$ as $n\rightarrow \infty$, we only place emphasis on the case where $y$ is small positive constant.
Bahadur-Rao  proved that, for given positive constant $y$,
\begin{eqnarray} \label{asdsfcsd}
\mathbf{P}\left( \frac{S_{n}}{n} > y \right)= \frac{  e^{- n \, \Lambda^\ast(y)} }{  \sigma_1\mbox{}_y   t_y   \sqrt{2\pi n }  } \left[ 1+ O(\frac{c_y }{n}) \right], \ \ \ \ \ n\rightarrow \infty,
\end{eqnarray}
where $c_y,$ $\sigma_1\mbox{}_y$ and $t_y$ depend on $y$ and the distribution of $\xi_1$; see also Bercu \cite{BR08,BR06}, Rozovky \cite{L05} and Gy\"{o}rfi,  Harrem\"{o}es and Tusn\'{a}dy \cite{GHT12} for more general results.
Our bound (\ref{fed}) implies that, for $y\geq 0$ small enough,
\begin{eqnarray}\label{adsdsvs}
\mathbf{P}\left( \frac{S_{n}}{n} > y \right)=   e^{   - n \, \Lambda^\ast(y)  }M(y\sqrt{n}) \Big[ 1+ O(y+\frac{1}{\sqrt{n}}) \Big].
\end{eqnarray}
In particular, when $0< y=y(n) \rightarrow 0 $ and $y \sqrt{n}\rightarrow \infty$ as $n\rightarrow \infty$,  we have
\begin{eqnarray} \label{adsgnsghs}
\mathbf{P}\left( \frac{S_{n}}{n} > y  \right)=  \frac{   e^{- n \, \Lambda^\ast(y )} }{   y    \sqrt{2\pi n }  } \Big[ 1+ o(1) \Big], \ \ \ \ \  \ n\rightarrow \infty.
\end{eqnarray}
Expansion (\ref{adsdsvs}) or (\ref{adsgnsghs}) is less precise than  (\ref{asdsfcsd}). However, the advantage of  the expansions (\ref{adsdsvs}) and (\ref{adsgnsghs}) over the Bahadur-Rao expansion (\ref{asdsfcsd}) is that  the expansions (\ref{adsdsvs}) or (\ref{adsgnsghs}) are uniform in $y$ (where $y$ may be dependent of $n$), in addition to the simpler expressions (without the factors $t_y$ and $\sigma_y$).

\section{Auxiliary results}\label{sec2.5}
We consider
the positive random variable
\[
Z_n(\lambda )=\prod_{i=1}^n\frac{e^{\lambda \xi _i}}{\mathbf{E}[ e^{\lambda \xi _i}] }%
,\ \ \ \ \ \ \ |\lambda|< \varepsilon^{-1},
\]
(the Esscher transformation) so that $\mathbf{E}[Z_{n}(\lambda)]=1$. We introduce the \emph{conjugate probability measure} $\mathbf{P}_\lambda $ defined by
\begin{equation}
d\mathbf{P}_\lambda =Z_n(\lambda )d\mathbf{P}.  \label{f21}
\end{equation}
Denote by $\mathbf{E}_{\lambda}$ the expectation with respect to $\mathbf{P}_{\lambda}$.
Setting
\[
b_i(\lambda)=\mathbf{E}_{\lambda}[ \xi_i ]=\frac{\mathbf{E} [\xi _ie^{\lambda \xi _i}] }{\mathbf{E}[e^{\lambda \xi _i}] },\quad i=1,...,n,
\]
and
\[
\eta_i(\lambda)=\xi_i - b_i(\lambda), \quad i=1, ...,n,
\]
we obtain the following decomposition:
\begin{equation}  \label{f22}
S_k=T_k(\lambda )+Y_k(\lambda ),\quad k=1, ...,n,
\end{equation}
where
\[
T_k(\lambda )=\sum_{i=1}^kb_i(\lambda ) \ \ \ \ \ \ \textrm{and} \ \ \ \ \ \ Y_k(\lambda )=\sum_{i=1}^k\eta _i(\lambda ).
\]

In the following, we give some lower and upper bounds of $T_n(\lambda ),$ which will be used in
the proofs of theorems.
\begin{lemma}
\label{lemma1} For all $0
\leq \lambda < \varepsilon^{-1} ,$
\begin{eqnarray*}
(1-2.4\lambda\varepsilon)\lambda\sigma^2 \leq  \frac{(1-1.5\lambda\varepsilon)(1-\lambda\varepsilon)}{1-\lambda\varepsilon+6
\lambda^2\varepsilon^2} \lambda\sigma^2  \leq T_n(\lambda ) \leq \frac{1
-0.5\lambda \varepsilon}{(1-\lambda\varepsilon)^2}\lambda\sigma^2.
\end{eqnarray*}
\end{lemma}

\noindent\emph{Proof.}  Since $\mathbf{E}[ \xi _i] =0$, by Jensen's inequality, we have
$\mathbf{E}[ e^{\lambda \xi _i}] \geq 1.$ Noting that
\[
\mathbf{E}[ \xi_{i} e^{\lambda\xi_{i}} ] =\mathbf{E} [\xi_{i}(e^{\lambda\xi_{i}}-1)] \geq 0,\ \ \ \ \lambda\geq 0,
\]
by Taylor's expansion of $e^x$, we get
\begin{eqnarray}
T_n(\lambda ) & \leq & \sum_{i=1}^{n}\mathbf{E} [\xi_{i} e^{\lambda \xi_{i}} ] \nonumber\\
& = & \lambda\sigma^2+ \sum_{i=1}^{n}\sum_{k=2}^{+\infty} \frac{\lambda^k}{k !}\mathbf{E} [\xi_{i}^{k+1} ]. \label{f24}
\end{eqnarray}
Using Bernstein's condition (\ref{Bernstein cond}), we obtain, for all $0 \leq \lambda < \varepsilon^{-1} ,$
\begin{eqnarray}
\sum_{i=1}^{n}\sum_{k=2}^{+\infty}\frac{\lambda^k}{k !}|\mathbf{E}[ \xi_{i}^{k+1} ]|  & \leq & \frac 12 \, \lambda^2 \sigma^2 \varepsilon
\sum_{k=2}^{+\infty}(k+1)(\lambda\varepsilon)^{k-2} \nonumber \\
& = & \frac{ 3-2\lambda\varepsilon }{2(1-\lambda\varepsilon)^2}\,
\lambda^2\sigma^2 \varepsilon .  \label{f25}
\end{eqnarray}
Combining (\ref{f24}) and (\ref{f25}), we get the desired upper bound of $T_{n}(\lambda)$.
By Jensen's inequality and \mbox{Bernstein's} condition (\ref{Bernstein cond}),
\begin{eqnarray*}
(\mathbf{E}[\xi_{i}^2 ])^{2} \leq \mathbf{E}[ \xi_{i}^4 ]\leq 12 \varepsilon^2 \mathbf{E}[\xi_{i}^2],
\end{eqnarray*}
from which we get
\[
\mathbf{E}[\xi_{i}^2] \leq 12 \varepsilon^2 .
\]
Using again Bernstein's condition (\ref{Bernstein cond}), we have, for all $0 \leq \lambda < \varepsilon^{-1},$
\begin{eqnarray}
\mathbf{E}[ e^{\lambda \xi_{i}} ]  & \leq & 1 + \sum_{k=2}^{+\infty}\frac{\lambda^k}{k !}|\mathbf{E} [\xi_{i}^k]| \nonumber\\
&\leq &1+ \frac{\lambda^2\mathbf{E}[\xi_{i}^2] }{2(1-\lambda\varepsilon)}  \nonumber \\
&\leq &1+ \frac{6\lambda^2\varepsilon^2 }{ 1-\lambda\varepsilon }  \nonumber \\
&= & \frac{1-\lambda\varepsilon+6\lambda^2\varepsilon^2 }{ 1-\lambda\varepsilon }.
\label{f26}
\end{eqnarray}
Notice that $g(t)=e^{t}-(1+t+\frac{1}{2}t^2)$ satisfies that $g(t)>0$ if $t>0$ and $g(t)<0$ if $t<0,$ which leads to  $tg(t)\geq0$ for all $t\in \mathbf{R}$. That is, $te^t \geq t(1+t+\frac12 t^2)$ for all $t\in \mathbf{R}$.
Therefore, for all $0 \leq \lambda < \varepsilon^{-1},$
$$\xi_ie^{\lambda\xi_i} \geq \xi_i \left(1+\lambda\xi_i+\frac{\lambda^2\xi_i^2}{2}\right).$$
Taking expectation, we get
\begin{eqnarray}
\mathbf{E}[\xi_{i} e^{\lambda \xi_{i}} ] \geq
\lambda \mathbf{E}[\xi_i^2]+ \frac{\lambda^{2}}{2}\mathbf{E}[\xi_{i}^{3}]
 \geq  \lambda \mathbf{E}[\xi_i^2] - \frac{\lambda^{2}}{2}\frac{1}{2}3!
\varepsilon \mathbf{E}[\xi_{i}^{2}]
 = (1-1.5 \lambda\varepsilon )\lambda \mathbf{E}[\xi_i^2],  \nonumber
\end{eqnarray}
from which, it follows that
\begin{eqnarray}
\sum_{i=1}^n\mathbf{E}[\xi_{i} e^{\lambda \xi_{i}}] &\geq& (1-1.5 \lambda\varepsilon )\lambda \sigma^2. \label{f27}
\end{eqnarray}
Combining (\ref{f26}) and (\ref{f27}), we obtain the following lower
bound of $T_n(\lambda )$: for all $0 \leq \lambda < \varepsilon^{-1},$
\begin{eqnarray}
T_n(\lambda )& \geq & \sum_{i=1}^{n} \frac{\mathbf{E}[\xi_{i} e^{\lambda \xi_{i}}]}{\mathbf{E}[e^{\lambda \xi_{i}}]}  \nonumber\\
&\geq& \frac{(1-1.5\lambda\varepsilon)(1-\lambda\varepsilon)}{1-\lambda\varepsilon+6
\lambda^2\varepsilon^2} \lambda\sigma^2  \nonumber \\
&\geq& (1-2.4\lambda\varepsilon)\lambda\sigma^2 . \label{f28}
\end{eqnarray}
This completes the proof of Lemma \ref{lemma1}.
\hfill\qed

We now consider the following \emph{cumulant} function
\begin{equation}
\Psi _n(\lambda )=\sum_{i=1}^n\log \mathbf{E}[e^{\lambda \xi _i}], \ \ \ \ \ \ \ \
0\leq \lambda < \varepsilon^{-1}.
\end{equation}
We have the following elementary bound for $\Psi_n
(\lambda ).$

\begin{lemma}
\label{lemma2} For all $0
\leq \lambda < \varepsilon^{-1} ,$
\begin{eqnarray}
\Psi _n(\lambda ) \leq n\log \left(1+ \frac{\lambda^2\sigma^2}{%
2n(1-\lambda\varepsilon)} \right) \leq \frac{\lambda^2\sigma^2}{2
(1-\lambda\varepsilon)}  \nonumber
\end{eqnarray}
and
\begin{eqnarray*}
-\lambda T_{n}(\lambda)+\Psi _n(\lambda ) &\geq& -\frac{\lambda^2\sigma^2}{%
2(1-\lambda\varepsilon)^6} .
\end{eqnarray*}
\end{lemma}

\noindent\emph{Proof.} By Bernstein's condition (\ref{Bernstein cond}), it is easy to see that, for all $0\leq \lambda < \varepsilon^{-1}$,
\begin{eqnarray*}
\mathbf{E}[ e^{\lambda \xi _i} ]  &=&1+ \sum_{k=2}^{+\infty}\frac{\lambda^{k}}{k !} \mathbf{E}
[\xi_{i}^{k}]
 \leq  1+ \frac{\lambda^2}{2}\mathbf{E} [\xi_{i}^{2}]
\sum_{k=2}^{\infty}(\lambda\varepsilon)^{k-2}
 =   1+ \frac{\lambda^2\mathbf{E} [\xi_{i}^{2}] }{2(1-\lambda\varepsilon)}.  \nonumber
\end{eqnarray*}
Then, we have
\begin{eqnarray}
\Psi _n(\lambda )  &\leq&\sum_{i=1}^{n} \log \left(1+ \frac{\lambda^2\mathbf{E}
[\xi_{i}^{2}]}{2(1-\lambda\varepsilon)} \right) . \label{f30}
\end{eqnarray}
Using the fact $\log(1+t)$ is concave in $t\geq 0$ and $\log (1+t) \leq t$, we get the first assertion of the lemma.
Since $\Psi _n(0)=0$ and $\Psi_n^{\prime }(\lambda)=T_{n}(\lambda)$, by Lemma \ref{lemma1}, for all $0 \leq
\lambda < \varepsilon^{-1} ,$
\begin{eqnarray}
\Psi _n(\lambda ) = \int_{0}^{\lambda}T_{n}(t)dt\geq
\int_{0}^{\lambda}t(1-2.4t\varepsilon)\sigma^2dt= \frac{\lambda^2\sigma^2}{2}(1-1.6\lambda \varepsilon)
 .  \nonumber
\end{eqnarray}
Therefore, using again Lemma \ref{lemma1}, we see that
\begin{eqnarray*}
-\lambda T_{n}(\lambda)+\Psi _n(\lambda ) &\geq& -\frac{1-0.5\lambda\varepsilon}{(1-\lambda\varepsilon)^2}\lambda^2 \sigma^2 + \frac{\lambda^2\sigma^2}{2}(1-1.6\lambda \varepsilon) \nonumber \\
&\geq& -\frac{\lambda^2\sigma^2}{2(1-\lambda\varepsilon)^6},
\end{eqnarray*}
which completes the proof of the second assertion of the lemma.
\hfill\qed

Denote $ \overline{\sigma}^2(\lambda)=\mathbf{E}_{\lambda}[Y_n^2(\lambda)]$.  By the relation between $\mathbf{E}$ and $\mathbf{E}_{\lambda}$, we have
\begin{eqnarray*}
\overline{\sigma}^2(\lambda)  =  \sum_{i=1}^n \left( \frac{\mathbf{E}[ \xi _i^2e^{\lambda \xi
_i}] }{\mathbf{E}[ e^{\lambda \xi _i} ]}-\frac{(\mathbf{E}[\xi _ie^{\lambda \xi _i} ])^2}{%
(\mathbf{E}[e^{\lambda \xi _i}] )^2}\right) , \ \ \ \ 0\leq \lambda < \varepsilon^{-1}.
\nonumber
\end{eqnarray*}

\begin{lemma}
\label{lemma5} For all $%
0\leq \lambda < \varepsilon^{-1},$
\begin{eqnarray}
 \frac{(1-\lambda\varepsilon)^2(1-3\lambda\varepsilon)}{(1-\lambda\varepsilon+6\lambda^2\varepsilon^2)^2 }\sigma^2   \leq  \overline{\sigma}^2(\lambda) \leq
\frac{\sigma^2}{(1-\lambda\varepsilon)^3}.
\end{eqnarray}
\end{lemma}

\noindent\emph{Proof.}
Denote $f(\lambda)=\mathbf{E}[ \xi_i^2e^{\lambda\xi_{i}}] \mathbf{E} [e^{\lambda \xi _i}] - (\mathbf{E}[\xi _ie^{\lambda \xi _i}])^2$.
Then,
\begin{eqnarray*}
f'(0) = \mathbf{E}[\xi_i^3]  \ \ \ \textrm{and} \ \ \
f''(\lambda)&=&\mathbf{E}[ \xi_i^4e^{\lambda\xi_{i}}] \mathbf{E}[ e^{\lambda \xi _i}] - (\mathbf{E}[\xi _i^2e^{\lambda \xi _i}])^2 \geq 0.
\end{eqnarray*}
Thus,
\begin{eqnarray}\label{f44}
f(\lambda) \geq f(0)+f'(0) \lambda = \mathbf{E}[\xi_i^2]+\lambda \mathbf{E}[\xi_{i}^3].
\end{eqnarray}
Using (\ref{f44}), (\ref{f26}) and Bernstein's condition (\ref{Bernstein cond}),
we have, for all $0 \leq \lambda < \varepsilon^{-1},$
\begin{eqnarray*}
 \mathbf{E}_{\lambda}[\eta_{i}^2 ]&=&  \frac{\mathbf{E}[ \xi _i^2e^{\lambda \xi _i} ]\mathbf{E}[ e^{\lambda \xi _i} ] -(\mathbf{E}[\xi _i e^{\lambda \xi _i}] )^2}{
(\mathbf{E}[e^{\lambda\xi _i}] )^2} \\
&\geq &  \frac{\mathbf{E}[\xi_i^2]+\lambda \mathbf{E}[\xi_{i}^3]}{
(\mathbf{E}[e^{\lambda\xi _i}] )^2}\\
&\geq & \left(\frac{1-\lambda\varepsilon}{1-\lambda\varepsilon+6\lambda^2\varepsilon^2 }\right)^2(\mathbf{E}[\xi_i^2]+\lambda \mathbf{E}[\xi_{i}^3]) \\
&\geq & \frac{(1-\lambda\varepsilon)^2(1-3\lambda\varepsilon)}{(1-\lambda\varepsilon+6\lambda^2\varepsilon^2)^2 }\mathbf{E}[\xi_i^2].
\end{eqnarray*}
Therefore
$$ \overline{\sigma}^2(\lambda)    \geq \frac{(1-\lambda\varepsilon)^2(1-3\lambda\varepsilon)}{(1-\lambda\varepsilon+6\lambda^2\varepsilon^2)^2 }\sigma^2. $$
Using Taylor's expansion of $e^x$ and Bernstein's condition (\ref{Bernstein cond}) again, we obtain $$\overline{%
\sigma}^2(\lambda)\leq \sum_{i=1}^n \mathbf{E}[ \xi _i^2e^{\lambda \xi _i} ] \leq
\frac{\sigma^2}{(1-\lambda\varepsilon)^3}.$$
This completes the proof of Lemma \ref{lemma5}.
\hfill\qed

For the random variable $Y_{n}(\lambda)$ with  $0\leq \lambda < \varepsilon^{-1}$, we have the following
result on the rate of convergence to the standard normal law.

\begin{lemma}
\label{lemma3} For all $%
0\leq \lambda < \varepsilon^{-1},$
\[
\sup_{y\in \mathbf{R}}\left| \mathbf{P}_\lambda \left( \frac{Y_n(\lambda )}{\overline{\sigma}(\lambda) }\leq
y \right)-\Phi (y)\right| \leq 13.44 \frac{\sigma^2 \varepsilon }{\overline{\sigma}%
^{3}(\lambda) (1-\lambda\varepsilon)^{4}} .
\]
\end{lemma}

\noindent\emph{Proof.} Since $Y_{n}(\lambda)=\sum_{i=1}^{n}\eta_{i}(\lambda)$ is the sum of independent and centered (respect to $\mathbf{P}_{\lambda}$) random variables $\eta_{i}(\lambda)$,
using standard results on  the rate of convergence in the central limit theorem (cf.\ e.g.\ Petrov \cite{Petrov75}, p.\ 115)  we get, for $0\leq \lambda < \varepsilon^{-1},$
\begin{eqnarray*}
 \sup_{y\in \mathbf{R}}\left| \mathbf{P}_\lambda \left(\frac{Y_n(\lambda )}{\overline{\sigma}(\lambda)} \leq y\right )-\Phi (y)\right| \leq C_1 \frac{1}{\overline{\sigma}^{3}(\lambda)}\sum_{i=1}^{n}\mathbf{E}_{\lambda} [|\eta_{i}|^{3}] ,
\end{eqnarray*}
where $C_1 >0$ is an absolute constant. For $0\leq \lambda < \varepsilon^{-1},$ using  Bernstein's condition, we have
\begin{eqnarray*}
\sum_{i=1}^{n}\mathbf{E}_{\lambda} [|\eta_{i}|^{3}]   &\leq &4\sum_{i=1}^{n}
\mathbf{E}_{\lambda}[|\xi_{i}|^3 + (\mathbf{E}_{\lambda}[|\xi_{i}|])^3 ]  \\
&\leq & 8 \sum_{i=1}^{n} \mathbf{E}_{\lambda} [|\xi_{i}|^3]   \\
&\leq & 8 \sum_{i=1}^{n} \mathbf{E} [|\xi_{i}|^3 \exp\{|\lambda\xi_{i}|\}]   \\
&\leq & 8 \sum_{i=1}^{n}\mathbf{E} \Big[ \sum_{j=0}^{\infty}\frac{\lambda^{j}}{j!}
|\xi_{i}|^{3+j} \Big]\\
&\leq & 4  \sigma^2\varepsilon \sum_{j=0}^{\infty}  (j+3)(j+2)(j+1)(\lambda\varepsilon)^j.
\end{eqnarray*}
As
\[
\sum_{j=0}^{\infty}  (j+3)(j+2)(j+1)x^j= \frac{d^3}{dx^3}\sum_{j=0}^{\infty} x^j=\frac{6}{(1-x)^4},\ \ \ |x|<1,
\]
we obtain, for  $0\leq \lambda < \varepsilon^{-1},$
\[
\sum_{i=1}^{n}\mathbf{E}_{\lambda} [|\eta_{i}|^{3} ]\leq 24 \frac{\sigma^2\varepsilon}{(1-\lambda\varepsilon)^4}.
\]
Therefore, we have, for $0\leq \lambda < \varepsilon^{-1},$
\begin{eqnarray*}
 \sup_{y\in \mathbf{R}}\left| \mathbf{P}_\lambda \left(\frac{Y_n(\lambda )}{\overline{\sigma}(\lambda)} \leq y\right )-\Phi (y)\right| &\leq& 24C_1 \ \frac{\sigma^2 \varepsilon }{\overline{\sigma}
^{3}(\lambda) (1-\lambda\varepsilon)^{4}} \\
&\leq& 13.44 \frac{\sigma^2 \varepsilon }{\overline{\sigma}
^{3}(\lambda) (1-\lambda\varepsilon)^{4}},
\end{eqnarray*}
where the last step holds as $C_1\leq 0.56$ (cf.\ Shevtsova \cite{S10}).
\hfill\qed

Using Lemma \ref{lemma3}, we  easily obtain the following lemma.

\begin{lemma}
\label{lemma4} For all $0\leq\lambda\leq 0.1\varepsilon^{-1}$,
\begin{eqnarray*}
\sup_{y\in \mathbf{R}}\left| \mathbf{P}_\lambda \left(Y_n(\lambda )\leq \frac{y\sigma}{
1-\lambda\varepsilon} \right)-\Phi (y)\right|
 \leq 1.07\lambda \varepsilon + 42.45\frac{\varepsilon}{\sigma}.
\end{eqnarray*}
\end{lemma}

\noindent\emph{Proof.}  Using Lemma \ref{lemma5}, we have, for all $0\leq \lambda < \frac13 \varepsilon^{-1}$,
\begin{eqnarray}
 \sqrt{1-\lambda\varepsilon} \ \leq\  \frac{ \sigma}{\overline{\sigma}(\lambda)(
1-\lambda\varepsilon)}\  \leq \ \frac{1-\lambda\varepsilon + 6\lambda^2\varepsilon^2}{(1-\lambda\varepsilon)^2\sqrt{1-3\lambda\varepsilon}}  . \label{fes}
\end{eqnarray}
It is easy to see that
\begin{eqnarray*}
 && \left| \mathbf{P}_\lambda \left(Y_n(\lambda )\leq \frac{y\sigma}{
1-\lambda\varepsilon} \right)-\Phi (y)\right|  \nonumber \\
&\leq& \left| \mathbf{P}_\lambda \left( \frac{Y_n(\lambda )}{\overline{\sigma}(\lambda)}\leq \frac{y\sigma}{\overline{\sigma}(\lambda)(
1-\lambda\varepsilon)} \right)-\Phi \left(\frac{y\sigma}{\overline{\sigma}(\lambda)(
1-\lambda\varepsilon)} \right)\right|  \\
&&+\left|  \Phi \left(\frac{y\sigma}{\overline{\sigma}(\lambda)(
1-\lambda\varepsilon)} \right) -\Phi (y)\right| \\
&=:& I_1+I_2.
\end{eqnarray*}
By Lemma \ref{lemma3} and (\ref{fes}), we get, for all $0\leq\lambda< \frac 13 \varepsilon^{-1}$,
\begin{eqnarray*}
I_1
\ \leq\ 13.44 \ \frac{\sigma^2 \varepsilon }{\overline{\sigma}
^{3}(\lambda) (1-\lambda\varepsilon)^{4}}\ \leq\ 13.44 R(\lambda\varepsilon) \frac{\varepsilon}{\sigma} .
\end{eqnarray*}
Using Taylor's expansion and (\ref{fes}), we obtain, for all $0\leq\lambda< \frac 13\varepsilon^{-1}$,
\begin{eqnarray*}
I_2 &\leq&  \frac{1}{\sqrt{2\pi}} y e^{-\frac{y^2(1-\lambda\varepsilon)}{2}}\left|\frac{\sigma}{\overline{\sigma}(\lambda)(1-\lambda\varepsilon)}-1\right|  \\
&\leq&  \frac{1}{\sqrt{2\pi}} y e^{-\frac{y^2(1-\lambda\varepsilon)}{2}} \left(\left|\frac{1-\lambda\varepsilon + 6\lambda^2\varepsilon^2}{(1-\lambda\varepsilon)^2\sqrt{1-3\lambda\varepsilon}}-1\right|\vee \left|1- \sqrt{1-\lambda\varepsilon} \right| \right)  \\
&\leq& \frac{1}{\sqrt{ 2 e \pi (1-\lambda\varepsilon)}}  \left|\frac{1-\lambda\varepsilon + 6\lambda^2\varepsilon^2}{(1-\lambda\varepsilon)^2\sqrt{1-3\lambda\varepsilon}}-1\right|.
\end{eqnarray*}
By simple calculations, we obtain,  for all $0\leq\lambda\leq0.1\varepsilon^{-1}$,
\begin{eqnarray*}
\left| \mathbf{P}_\lambda \left(Y_n(\lambda )\leq \frac{y\sigma}{
1-\lambda\varepsilon} \right)-\Phi (y)\right|
 \leq 1.07\lambda \varepsilon + 42.45\frac{\varepsilon}{\sigma}.
\end{eqnarray*}
This completes the proof of Lemma \ref{lemma4}.
\hfill\qed

\section{Proofs of Theorems \ref{th2}-\ref{th3}}\label{sec3.5}
In this section, we give upper bounds for $\mathbf{P}(S_n>x\sigma)$. For all $x\geq0$ and $0\leq \lambda <\varepsilon^{-1},$
by (\ref{f21}) and (\ref{f22}), we have:
\begin{eqnarray}
\mathbf{P}(S_n>x\sigma )
 &= & \mathbf{E}_\lambda [Z_n (\lambda)^{-1}\mathbf{1}_{\{S_n>x\sigma \}}] \nonumber\\
&= & \mathbf{E}_\lambda [ e^{ -\lambda S_n+\Psi _n(\lambda ) } \mathbf{1}_{\{S_n>x\sigma \}} ]\nonumber \\
&= & \mathbf{E}_\lambda [e^{-\lambda T_{n}(\lambda) +\Psi _n(\lambda
)-\lambda Y_n(\lambda) } \mathbf{1}_{\{ Y_{n}(\lambda)+T_n(\lambda)-x\sigma >0 \}} ]\label{gnk}.
\end{eqnarray}
Setting $U_{n}(\lambda)=\lambda (Y_{n}(\lambda)+T_{n}(\lambda)-x\sigma)$, we get
\begin{eqnarray}
\mathbf{P}(S_n>x\sigma ) &=& e^{-\lambda x\sigma +\Psi _n(\lambda )
 } \mathbf{E}_{\lambda} [e^{- U_{n}(\lambda)}\mathbf{1}_{\{ U_{n}(\lambda)> 0\}}].\nonumber
\end{eqnarray}
Then, we deduce, for all $x\geq 0$ and $0\leq \lambda <\varepsilon^{-1},$
\begin{eqnarray}
\mathbf{P}(S_n>x\sigma ) &=& e^{- \lambda x\sigma +\Psi _n(\lambda ) } \int_0^\infty e^{-t} \mathbf{P}_{ \lambda}
(0< U_{n}( \lambda ) \leq t) dt.\label{f46}
\end{eqnarray}
In the sequel, denote by $N(0, 1)$ a standard normal random variable.

\subsection{Proof of Theorem \ref{th2}} From (\ref{f46}), using Lemma \ref{lemma2}, we obtain, for all $x\geq0$ and $0\leq \lambda < \varepsilon^{-1} ,$
\begin{eqnarray}
 \mathbf{P}(S_n>x\sigma )
 &\leq & e^{-\lambda x\sigma + \frac{\lambda^2\sigma^2}{2
(1-\lambda\varepsilon)}  }  \int_0^\infty e^{-t} \mathbf{P}_{ \lambda}
(0< U_{n}( \lambda ) \leq t) dt.\label{fg33}
\end{eqnarray}
For any $x\geq0$ and $\beta \in [0, 0.5)$, let $\overline{\lambda }=\overline{\lambda }(x)\in [0,\varepsilon^{-1})$ be the
unique solution of the equation
\[
\frac{\lambda-\beta\lambda^2\varepsilon }{(1-\lambda\varepsilon)^2}=\frac{x}{\sigma }.
\]
This definition and Lemma \ref{lemma1} implies that
\begin{equation}
\overline{\lambda} =\frac{2x/\sigma }{1+2x\varepsilon/\sigma+\sqrt{1 +4(1- \beta)x \varepsilon /\sigma }%
} \ \ \ \ \ \mbox{and} \ \ \ \   T_{n}(
\overline{\lambda}) \leq x\sigma. \label{kj}
\end{equation}
Using (\ref{fg33}) with $\lambda =\overline{\lambda }$,  we get
\begin{eqnarray}
 \ \ \ \ \ \ \mathbf{P}(S_n>x\sigma )
 \leq  e^{- \frac{1}{2 }\left(1+ (1-2\beta)\overline{\lambda}\varepsilon \right)\,\widetilde{x} \,^2  } \int_0^\infty e^{-t} \mathbf{P}_{\overline{\lambda}}(0< U_{n}(\overline{\lambda}) \leq t) dt, \label{f35s}
\end{eqnarray}
where $$\widetilde{x}=\frac{\overline{\lambda}\sigma }{1-\overline{\lambda}\varepsilon}.$$
By (\ref{kj}) and Lemma {\ref{lemma4}}, we have, for $0\leq \overline{\lambda}\leq 0.1\varepsilon^{-1} $,
\begin{eqnarray}
&&\int_0^\infty e^{-t} \mathbf{P}_{\overline{\lambda}}(0< U_{n}(\overline{\lambda}) \leq t) dt \nonumber\\
&=&\int_0^\infty e^{ -y\widetilde{x}  }\mathbf{P}_{\overline{\lambda}}\left(0< U_{n}(\overline{\lambda}) \leq y\widetilde{x}  \right)  \widetilde{x}  dy  \nonumber \\
&\leq &\int_0^\infty e^{-y\widetilde{x} } \mathbf{P} \left(0< N(0,1) \leq y \right) \widetilde{x} dy   +  2 \left( 1.07\overline{\lambda} \varepsilon + 42.45 \frac{\varepsilon}{\sigma}   \right)   \nonumber \\
&=& M\left(\widetilde{x}\right) +   2.14 \overline{\lambda} \varepsilon + 84.9 \frac{\varepsilon}{\sigma}. \label{f36}
\end{eqnarray}
Since $\int_0^\infty e^{-t} \mathbf{P}_{\overline{\lambda}}(0< U_{n}(\overline{\lambda}) \leq t) dt\leq 1$ and $ M^{-1}(t)\leq  \sqrt{2\pi} \left( 1 + t\right)$ for $t\geq0$ (cf.\ (\ref{mt})),
combining (\ref{f35s}) and (\ref{f36}), we deduce, for all $x\geq 0 ,$
\begin{eqnarray}
\ \ && \mathbf{P}(S_n>x\sigma ) \nonumber\\
&\leq& e^{-\frac{1}{2}(1-2\beta)\overline{\lambda}\varepsilon \widetilde{x}\,^2 - \frac{1}{2}\widetilde{x}\,^2 }\mathbf{1}_{\{\overline{\lambda}\varepsilon> 0.1\}}  \nonumber\\
&& +e^{-\frac{1}{2} (1-2\beta)\overline{\lambda}\varepsilon  \widetilde{x}\,^2  }\left[ 1-\Phi\left(\widetilde{x} \right) +
 e^{-\frac{1}{2}\widetilde{x}\,^2  } \left(2.14\overline{\lambda}\varepsilon +84.9\frac{\varepsilon}{\sigma}\right)\right] \mathbf{1}_{\{ \overline{\lambda}\varepsilon \leq 0.1\}}  \nonumber\\
&\leq&\left(1-\Phi\left(\widetilde{x} \right)\right) ( I_{11}+ I_{12} ), \label{fmou}
\end{eqnarray}
with
\begin{eqnarray}
 I_{11}= \exp\left\{-\frac{1 }{2 }(1-2\beta)\overline{\lambda}\varepsilon\widetilde{x}^2\right\}  \left[ \sqrt{2\pi}\left(1+\widetilde{x} \right)
  \right]\mathbf{1}_{\{\overline{\lambda}\varepsilon> 0.1\}}
\end{eqnarray}
and
\begin{eqnarray*}
  I_{12}&=&  e^{-\frac{1}{2 }(1-2\beta)\overline{\lambda}\varepsilon\widetilde{x}\,^2 }  \left[ 1+ \sqrt{2\pi}\left(1+\widetilde{x} \right) \left( 2.14 \overline{\lambda} \varepsilon + 84.9 \frac{\varepsilon}{\sigma}  \right)  \right]\mathbf{1}_{\{ \overline{\lambda}\varepsilon\leq 0.1\}}.
\end{eqnarray*}
Now we shall give estimates for $I_{11}$ and $I_{12}$.
If $ \overline{\lambda}\varepsilon > 0.1$,
then $I_{12}=0$ and
\begin{eqnarray}
 I_{11} &\leq& \exp\left\{-0.1(1-2\beta)\frac{ \widetilde{x}^2}{2 }\right\}  \left[ \sqrt{2\pi}\left(1+\widetilde{x} \right)
  \right]  .
\end{eqnarray}
By a simple calculation, $I_{11} \leq 1$ provided that $\widetilde{x}\geq \frac{8}{1-2\beta}$ (note that $\beta \in [0, 0.5)$). For $0\leq \widetilde{x}< \frac{8}{1-2\beta}$, we get
$\overline{\lambda}\sigma = \widetilde{x}(1- \overline{\lambda} \varepsilon) < \frac{8}{1-2\beta} (1-0.1)=\frac{7.2}{1-2\beta}$. Then, using $10\overline{\lambda}\varepsilon>1$, we obtain
\begin{eqnarray*}
 I_{11} &\leq& 1+ \sqrt{2\pi}\left(1+\widetilde{x} \right) \\
 &\leq& 1+ 10\sqrt{2\pi}\left(1+\widetilde{x} \right)\overline{\lambda}\sigma \frac \varepsilon  \sigma\\
 &\leq& 1+  \frac{180.48}{1-2\beta} \left(1+\widetilde{x} \right) \frac \varepsilon  \sigma.
\end{eqnarray*}
If $0\leq \overline{\lambda}\varepsilon \leq 0.1$, we have $I_{11}=0$. Since
\begin{eqnarray*}
&&1+ \sqrt{2\pi}\left(1+\widetilde{x} \right) \left( 2.14 \overline{\lambda} \varepsilon + 84.9 \frac{\varepsilon}{\sigma}   \right)  \\
&\leq& \left(1+  2.14\sqrt{2\pi}\left(1+\widetilde{x} \right)  \overline{\lambda} \varepsilon \right)\left(1+ 84.9\sqrt{2\pi}\left(1+\widetilde{x} \right)  \frac{\varepsilon}{\sigma}  \right)\\
&=& J_1  J_2,
\end{eqnarray*}
it follows that $I_{12}  \leq  \exp\left\{-\frac{1}{2 }(1-2\beta)\overline{\lambda}\varepsilon \widetilde{x}^2\right\} J_1   J_2$.
Using the inequality $1+x\leq e^x$, we deduce
\begin{eqnarray*}
I_{12} &\leq& \exp\left\{-\overline{\lambda}\varepsilon  \left( (1-2\beta) \frac{ \widetilde{x}^2}{2 }- 2.14\sqrt{2\pi}\left(1+\widetilde{x} \right) \right)\right\}J_2.
\end{eqnarray*}
If $\widetilde{x}\geq \frac{11.65}{1-2\beta}$, we see that $ \frac{1}{2 }(1-2\beta)\widetilde{x}^2- 2.14\sqrt{2\pi}\left(1+\widetilde{x} \right) \geq 0$, so $I_{12}\leq J_2$.
For $0\leq \widetilde{x}< \frac{11.65}{1-2\beta}$, we get
$\overline{\lambda}\sigma = \widetilde{x}(1- \overline{\lambda} \varepsilon) < \frac{11.65}{1-2\beta}$. Then
\begin{eqnarray*}
 I_{12} &\leq& 1+ \sqrt{2\pi}\left(1+\widetilde{x} \right) \left( 2.14 \overline{\lambda} \varepsilon + 84.9 \frac{\varepsilon}{\sigma}   \right) \\
  &<& 1+ \sqrt{2\pi}\left(1+\widetilde{x} \right) \left( 2.14 \frac{11.65}{1-2\beta} + 84.9  \right) \frac \varepsilon  \sigma\\
 &\leq& 1+ \left(  \frac{62.493}{1-2\beta} + 212.813  \right) \left(1+\widetilde{x} \right) \frac \varepsilon  \sigma.
\end{eqnarray*}
Hence, whenever $0\leq \overline{\lambda}\varepsilon < 1$, we have
\begin{eqnarray}
I_{11}+I_{12} &\leq&1+  \left( \left(  \frac{62.493}{1-2\beta} + 212.813  \right) \vee \frac{180.48}{1-2\beta}  \right)\left(1+\widetilde{x} \right)   \frac{\varepsilon}{\sigma} . \label{i1i2}
\end{eqnarray}
Therefore, substituting $\overline{\lambda}$ from (\ref{kj}) in the expression of $\widetilde{x}=\frac{\overline{\lambda}\sigma}{1-\overline{\lambda}\varepsilon}$ and replacing $1-2\beta$ by $\delta$, we obtain  inequality (\ref{f10}) in Theorem \ref{th2} from (\ref{fmou}) and (\ref{i1i2}).

\subsection{Proof of Theorem \ref{th3}} For any $x\geq 0$,
let $\overline{\lambda }=\overline{\lambda }(x)\in [0,\varepsilon^{-1})$ be the
unique solution of the equation
\begin{equation}
\frac{\lambda-0.5\lambda^2\varepsilon}{(1-\lambda\varepsilon)^2}=\frac{x}{\sigma}.
\end{equation}
By Lemma \ref{lemma1}, it follows that
\begin{equation}
\overline{\lambda} =\frac{2x/\sigma }{1+2x\varepsilon/\sigma+\sqrt{1+2x\varepsilon/\sigma }
 }\ \ \ \ \mbox{and} \ \ \ \ \ T_{n}(\overline{
\lambda})\leq x\sigma.  \label{fkjl}
\end{equation}
  Using Lemma {\ref{lemma3}} and $T_{n}(\overline{
\lambda})\leq x\sigma$, we have,   for all $0\leq \overline{\lambda} < \varepsilon^{-1} ,$
\begin{eqnarray}
&& \int_0^\infty e^{-t} \mathbf{P}_{ \overline{\lambda}}
(0< U_{n}( \overline{\lambda} ) \leq t) dt \nonumber\\
&=&\int_0^\infty e^{-y \overline{\lambda}\overline{\sigma}(\overline{\lambda})}\mathbf{P}_{\overline{\lambda}}\left(0<U_{n}(\overline{\lambda}) \leq y \overline{\lambda} \overline{\sigma}(\overline{\lambda}) \right) \overline{\lambda} \overline{\sigma}(\overline{\lambda})dy \nonumber\\
&\leq&\int_0^\infty e^{-y \overline{\lambda}\overline{\sigma}(\overline{\lambda})} \mathbf{P} \left(0< N(0, 1)\leq y  \right) \overline{\lambda} \overline{\sigma}(\overline{\lambda})dy+ 26.88 \, \frac{\sigma^2 \varepsilon }{\overline{\sigma}^{3}(\overline{\lambda})
(1-\overline{\lambda}\varepsilon)^{4}} \nonumber\\
&\leq &\int_0^\infty e^{- y\overline{\lambda} \overline{\sigma}(\overline{\lambda})}d \Phi(y)   + 26.88 \ \frac{\sigma^2
\varepsilon }{\overline{\sigma}^{3}(\overline{\lambda}) (1-\overline{\lambda}\varepsilon)^{4}}  \nonumber \\
&=&\ F\ :=\  M\left(\overline{\lambda} \overline{\sigma}(\overline{\lambda}) \right) + 26.88 \ \frac{\sigma^2
\varepsilon }{\overline{\sigma}^{3}(\overline{\lambda}) (1-\overline{\lambda}\varepsilon)^{4}}. \label{sfka}
\end{eqnarray}
 Using $\lambda=\overline{\lambda}$ and $\int_0^\infty e^{-t} \mathbf{P}_{ \overline{\lambda}}
(0< U_{n}( \overline{\lambda} ) \leq t) dt\leq 1$, from  (\ref{f46}) and (\ref{sfka}), we obtain
\begin{eqnarray} \label{fkf}
\ \ \ \ \ \mathbf{P}(S_n>x\sigma ) &\leq& \left[F\wedge 1 \right] \times \exp \left\{- \overline{\lambda} x\sigma +\Psi _n(\overline{\lambda} )\right\}
    .
\nonumber
\end{eqnarray}
By Lemma \ref{lemma2}, inequality (\ref{fkf}) implies that
\begin{eqnarray}
\mathbf{P}(S_n>x\sigma ) &\leq& \left[F\wedge 1 \right] \times \exp \left\{- \overline{\lambda} x\sigma + n\log \left(1+ \frac{%
\overline{\lambda}^2\sigma^2}{2n(1-\overline{\lambda}\varepsilon)} \right)
\right\}  .  \nonumber
\end{eqnarray}
Substituting $\overline{\lambda}$ from (\ref{fkjl})  in the previous exponential function, we get
\begin{eqnarray}
  \mathbf{P}(S_n> x\sigma )
   &\leq& \left[F\wedge 1 \right] \times B_{n}\left(x, \frac{\varepsilon}\sigma \right) .   \label{fhk}
\end{eqnarray}
Next, we give an estimation of $F$.
Since $M(t)$ is decreasing
in $t\geq0$ and $|M'(t)|\leq  \frac{1}{\sqrt{\pi}\, t^2}, t>0$, it follows that
\begin{eqnarray}
M\left(\overline{\lambda} \overline{\sigma}(\overline{\lambda}) \right) -M(x)  \nonumber
&\leq& \frac{1}{\sqrt{\pi}}\frac{1}{\overline{\lambda}^2 \overline{\sigma}^2(\overline{\lambda}) } \left( x - \overline{\lambda} \overline{\sigma}(\overline{\lambda})\right)^+ . \nonumber
\end{eqnarray}
Using Lemma \ref{lemma5}, by a simple calculation, we  deduce
\begin{eqnarray}
 && M\left(\overline{\lambda} \overline{\sigma}(\overline{\lambda}) \right) -M(x) \nonumber\\
 &\leq& \frac{1}{\sqrt{\pi}} \frac{ \overline{\lambda}\sigma }{\overline{\lambda}^2 \overline{\sigma}^2(\overline{\lambda})} \left( \frac{1-0.5\overline{\lambda}\varepsilon}{(1-\overline{\lambda}\varepsilon)^2} -    \frac{(1-\overline{\lambda}\varepsilon)\sqrt{ (1-3\overline{\lambda}\varepsilon)^+  }}{1-\overline{\lambda}\varepsilon+6\overline{\lambda}^2\varepsilon^2} \right)  \nonumber \\
&\leq&  \frac{(1-0.5\overline{\lambda} \varepsilon)(1-\overline{\lambda}\varepsilon+6\overline{\lambda}^2\varepsilon^2)-(1-\overline{\lambda}\varepsilon) \sqrt{(1-3\overline{\lambda}\varepsilon)^{+}}  }{\sqrt{\pi} \ \overline{\lambda} \varepsilon (1-\overline{\lambda}\varepsilon)^4(1-3\overline{\lambda}\varepsilon)^+/(1-\overline{\lambda}\varepsilon +6\overline{\lambda}^2\varepsilon^2)} \, \frac{\varepsilon}{\sigma }  \nonumber\\
&\leq& 1.11 R\left(\overline{\lambda}\varepsilon \right)  \frac{\varepsilon}{\sigma}. \label{fkml}
\end{eqnarray}
By Lemma \ref{lemma5}, it is easy to see that
\begin{equation}\label{fnkmavsf}
26.88 \,\frac{\sigma^2 \varepsilon }{\overline{%
\sigma}^{3}(\overline{\lambda}) (1- \overline{\lambda} \varepsilon)^{4}} \leq 26.88 R\left(\overline{\lambda}\varepsilon \right) \frac{\varepsilon}{\sigma}.
\end{equation}
Hence, it follows from (\ref{sfka}), (\ref{fkml}) and (\ref{fnkmavsf}) that
\begin{eqnarray}\label{fsj}
 F \ \leq \ M(x)+ 27.99 R\left( \overline{\lambda} \varepsilon \right) \frac{\varepsilon}{\sigma}.
\end{eqnarray}
Implementing (\ref{fsj}) into (\ref{fhk}) and using $\overline{\lambda}\varepsilon \leq x\frac{\varepsilon}{\sigma}$, we obtain inequality (\ref{f17}).

\section{Proof of Theorem  \ref{th4}}\label{sec3}
In this section, we give a
lower bound for $\mathbf{P}(S_n>x\sigma)$.
From Lemma \ref{lemma2} and (\ref{gnk}), it follows that, for all $0\leq \lambda < \varepsilon^{-1},$
\begin{eqnarray}  \label{f48}
\mathbf{P}(S_n>x\sigma ) &\geq &\exp\left\{-\frac{\lambda^2\sigma^2}{%
2(1-\lambda\varepsilon)^6} \right\} \mathbf{E}_\lambda [ e^{-\lambda Y_n(\lambda)}
\mathbf{1}_{\{ Y_{n}(\lambda)+T_n(\lambda)-x\sigma >0\}}].
\nonumber
\end{eqnarray}
Let $\overline{\lambda }=\overline{\lambda }(x)\in [0,\varepsilon^{-1}/4.8]$ be
the unique solution of the equation
\begin{equation}
\lambda(1-2.4\lambda\varepsilon)\sigma^2 =x\sigma .
\end{equation}
This definition and Lemma \ref{lemma1} implies that, for all $0\leq x \leq
\sigma/(9.6\varepsilon)$,
\begin{equation}
\overline{\lambda} =\frac{2x/\sigma }{1 +\sqrt{1-9.6 x \varepsilon /\sigma }} \
\ \ \ \ \mbox{and} \ \ \ \ x\sigma \leq T_{n}(\overline{\lambda}).
\end{equation}
Therefore,
\begin{eqnarray}
\mathbf{P}(S_n>x\sigma ) &\geq &\exp\left\{-\frac{\lambda^2\sigma^2}{%
2(1-\lambda\varepsilon)^6} \right\} \mathbf{E}_\lambda [e^{-\lambda Y_n(\lambda)}
\mathbf{1}_{\{ Y_{n}(\lambda)  >0\}}].
\nonumber
\end{eqnarray}
Setting $V_{n}(\overline{\lambda})=\overline{\lambda} Y_{n}(\overline{\lambda}),$ we get
\begin{eqnarray}  \label{f50}
\mathbf{P}(S_n>x\sigma ) &\geq& \exp \left\{-\frac{\check{x}^2}{2 } \right\} \int_0^\infty e^{-t} \mathbf{P}_{ \overline{\lambda}}
(0< V_{n}( \overline{\lambda} ) \leq t) dt,
\end{eqnarray}
where $\check{x}=\frac{\overline{\lambda} \sigma}{(1-\overline{\lambda}\varepsilon)^3}$.
By Lemma {\ref{lemma3}} and an argument similar to that used to prove  (\ref{sfka}), it is easy to see that
\begin{eqnarray*}
  \int_0^\infty e^{-t} \mathbf{P}_{ \overline{\lambda}}
(0< V_{n}( \overline{\lambda} ) \leq t) dt
&\geq& M\left( \overline{\lambda}\overline{\sigma}(\overline{\lambda}) \right) - G ,  \nonumber
\end{eqnarray*}
where $G=26.88 \ \frac{\sigma^2 \varepsilon }{\overline{\sigma}%
^{3}(\overline{\lambda}) (1-\overline{\lambda}\varepsilon)^{4}}. $
Since $M(t)$ is decreasing in $t\geq0$ and $\overline{\sigma}(\overline{\lambda})
\leq \frac{\sigma}{(1-\overline{\lambda}\varepsilon)^3}$ (cf.\ Lemma \ref{lemma5}), it follows that
\begin{eqnarray}
\int_0^\infty e^{-t} \mathbf{P}_{ \overline{\lambda}}
(0< V_{n}( \overline{\lambda} ) \leq t) dt &\geq& M\left( \check{x}\right)   - G .  \nonumber
\end{eqnarray}
Returning to (\ref{f50}), we obtain
\begin{eqnarray}
\mathbf{P}(S_n>x\sigma ) &\geq& 1-\Phi\left(\check{x} \right)  -    G \exp \left\{- \frac{\check{x}^2}{2}\right\}
.  \nonumber
\end{eqnarray}
Using Lemma \ref{lemma5}, for all $0\leq x \leq \sigma/(9.6\varepsilon)$, we have $0\leq \overline{\lambda}\varepsilon \leq 1/4.8$ and
\[
G  \geq 26.88 R\left( \overline{\lambda}\varepsilon\right)\frac{\varepsilon }{\sigma}.
\]
Therefore, for all $0\leq x \leq
\sigma/(9.6\varepsilon),$
\begin{eqnarray}
\ \  \mathbf{P}(S_n>x\sigma ) &\geq&   1-\Phi\left(\check{x}\right)  -26.88 R\left( \overline{\lambda}\varepsilon\right)\frac{\varepsilon }{\sigma}  \exp \left\{- \frac{\check{x}^2}{2} \right\}.  \nonumber
\end{eqnarray}
Using the inequality $M^{-1} (t) \leq \sqrt{2\pi} \left( 1+t \right)$ for $t\geq 0,$
we get, for all $0\leq x \leq \sigma/(9.6\varepsilon)$,
\begin{eqnarray*}
\mathbf{P}(S_n>x\sigma ) &\geq& \Big( 1-\Phi\left(\check{x} \right) \Big) \left[1 - 67.38R\left( \overline{\lambda}\varepsilon\right)\left(1+\check{x} \right)\frac{\varepsilon }{\sigma}\right].
\end{eqnarray*}
In particular, for all $0\leq x \leq \alpha \sigma/ \varepsilon$ with $0\leq \alpha\leq 1/9.6$, a simple calculation shows that
\[
0\leq \overline{\lambda}\varepsilon \leq \frac{2\alpha}{1+\sqrt{1-9.6 \alpha}}\leq \frac{1}{4.8}
\]
and
\[
67.38R\left( \overline{\lambda}\varepsilon\right)\leq 67.38R\left( \frac{2\alpha}{1+\sqrt{1-9.6 \alpha}} \right) \leq 67.38R\left( \frac{1}{4.8}\right) \leq 1753.23.
\]
This completes the proof of Theorem \ref{th4}.

\section{Proof of Theorem \ref{thend}}\label{secend}
Notice that $\Psi_n'(\lambda ) = T_{n}( \lambda)\in [0, \infty)$ is nonnegative in $\lambda\geq0$. Let $\overline{\lambda }=\overline{\lambda }(x)\geq0$ be the
unique solution of the equation $x\sigma = \Psi_n'(\lambda )$.
This definition implies that $T_{n}(\overline{\lambda})=   x\sigma$, $U_{n}(\overline{\lambda})=\overline{\lambda} Y_{n}(\overline{\lambda})$ and
\begin{eqnarray}
e^{- \overline{\lambda} x\sigma +\Psi _n(\overline{\lambda} ) } \ \ =\ \ \inf_{\lambda\geq0} e^{- \lambda  x\sigma +\Psi _n( \lambda  ) }
\ \ = \ \  \inf_{\lambda\geq0}\mathbf{E}[e^{\lambda(S_n-x\sigma)}].
\end{eqnarray}
From (\ref{f46}), using Lemma {\ref{lemma3}} with $\lambda =\overline{\lambda }$ and an argument similar to (\ref{sfka}),  we obtain
\begin{eqnarray}   \label{fkns}
\ \ \ \ \ \mathbf{P}(S_n>x\sigma ) &=&  \left( M\left(\overline{\lambda} \overline{\sigma}(\overline{\lambda}) \right)+
 \ \frac{26.88\theta_1 \sigma^2 \varepsilon }{\overline{\sigma}^{3}(\overline{\lambda})
(1-\overline{\lambda}\varepsilon)^{4}} \right) \inf_{\lambda\geq 0} \textbf{E}[e^{\lambda(S_n-x\sigma)}],\ \ \ \
\end{eqnarray}
where $|\theta_1|\leq 1$. Since $M(t)$ is decreasing
in $t\geq0$ and $|M'(t)|\leq  \frac{1}{\sqrt{\pi}\, t^2}$ in $ t>0$, it follows that
\begin{eqnarray}
\left| M\left(\overline{\lambda} \overline{\sigma}(\overline{\lambda}) \right) -M(x)  \right|
&\leq& \frac{1}{\sqrt{\pi}}\frac{ \left| x - \overline{\lambda} \overline{\sigma}(\overline{\lambda})\right|}{\overline{\lambda}^2 \overline{\sigma}^2(\overline{\lambda})\wedge x^2 }  . \label{jks}
\end{eqnarray}
By Lemma \ref{lemma1}, we have the following two-sided bound of $x$:
\begin{eqnarray}
 \frac{(1-1.5\overline{\lambda}\varepsilon)(1-\overline{\lambda}\varepsilon)}{1-\overline{\lambda}\varepsilon+6
\overline{\lambda}^2\varepsilon^2} \overline{\lambda}\sigma   \leq \frac{T_{n}(\overline{\lambda})}{\sigma}=x \leq \frac{1-0.5\overline{\lambda}\varepsilon}{(1-\overline{\lambda}\varepsilon)^2}\, \overline{\lambda}\sigma . \label{tbound}
\end{eqnarray}
Using the two-sided bound in Lemma \ref{lemma5} and (\ref{tbound}), by a simple calculation, we  deduce
\begin{eqnarray}\label{fxram}
 \overline{\lambda}^2 \overline{\sigma}^2(\overline{\lambda})\wedge x^2 \geq \frac{(1-\overline{\lambda}\varepsilon)^2(1-3\overline{\lambda}\varepsilon)}{(1-\overline{\lambda}\varepsilon+6\overline{\lambda}^2\varepsilon^2)^2}\overline{\lambda}^2\sigma^2
\end{eqnarray}
and
\begin{eqnarray}
 \left| x - \overline{\lambda} \overline{\sigma}(\overline{\lambda})\right| \leq  \overline{\lambda}\sigma \left( \frac{1-0.5\overline{\lambda}\varepsilon}{(1-\overline{\lambda}\varepsilon)^2} -   \frac{(1-\overline{\lambda}\varepsilon)\sqrt{ (1-3\overline{\lambda}\varepsilon)^+  }}{1-\overline{\lambda}\varepsilon+6\overline{\lambda}^2\varepsilon^2} \right).\label{sfgd}
 \end{eqnarray}
From (\ref{jks}), (\ref{fxram}), (\ref{sfgd}) and Lemma \ref{lemma5},  we easily obtain
\begin{eqnarray}\label{fbi2}
\left| M\left(\overline{\lambda} \overline{\sigma}(\overline{\lambda}) \right) -M(x) \right|
&\leq& 1.11  R\left(\overline{\lambda}\varepsilon \right)  \frac{\varepsilon}{\sigma}.
\end{eqnarray}
By Lemma \ref{lemma5}, it is easy to see that
\begin{equation}\label{fbi1}
 \frac{26.88 \sigma^2 \varepsilon }{\overline{%
\sigma}^{3}(\overline{\lambda}) (1- \overline{\lambda} \varepsilon)^{4}} \leq 26.88  R\left(\overline{\lambda}\varepsilon \right) \frac{\varepsilon}{\sigma}.
\end{equation}
Combining (\ref{fbi2}) and (\ref{fbi1}), we get, for all $0\leq \overline{\lambda } < \frac13 \varepsilon^{-1}$,
\begin{eqnarray}\label{fedsj}
  M\left(\overline{\lambda} \overline{\sigma}(\overline{\lambda}) \right) +  \frac{26.88 \theta_1 \sigma^2 \varepsilon }{\overline{%
\sigma}^{3}(\overline{\lambda}) (1- \overline{\lambda} \varepsilon)^{4}}
= M(x)+ 27.99 \theta_2 R\left( \overline{\lambda} \varepsilon \right) \frac{\varepsilon}{\sigma},
\end{eqnarray}
where $|\theta_2|\leq 1$.
By (\ref{tbound}), it follows that, for all $0\leq \overline{\lambda } < \frac13 \varepsilon^{-1}$,
\begin{eqnarray}\label{fedddfsj}
\overline{\lambda}\varepsilon \leq   \frac{1-\overline{\lambda}\varepsilon+6
\overline{\lambda}^2\varepsilon^2}{(1-1.5\overline{\lambda}\varepsilon)(1-\overline{\lambda}\varepsilon)}  x\frac{\varepsilon}{\sigma} \leq 4x\frac{\varepsilon}{\sigma}.
\end{eqnarray}
Implementing (\ref{fedsj}) into (\ref{fkns}) and using  (\ref{fedddfsj}), we obtain equality (\ref{fed}) of Theorem \ref{thend}.
Notice that $R < \infty$ restricts  $0\leq4x\frac{\varepsilon}{\sigma} < \frac1 3,$ which implies that $ 0\leq x < \frac1{12}\frac{\varepsilon}{\sigma}.$

\Acknowledgements{We would like to thank the referees for their helpful comments and suggestions.
Fan was partially supported by the Post-graduate Study Abroad Program sponsored by China Scholarship Council.
Liu was partially supported by the National Natural Science Foundation of China (Grant no.\ 11171044 and no.\ 11401590).
}



\begin{thebibliography}{99}
\bahao\baselineskip 11.5pt
\bibitem{A89} Arkhangelskii A N. Lower bounds for probabilities of large deviations for sums of independent random variables. Theory Probab Appl, 1989,  34: 565--575

\bibitem{BR60} Bahadur R, Rao R R.  On deviations of the sample mean. Ann Math Statist, 1960, 31: 1015--1027


\bibitem{Be62} Bennett G. Probability inequalities for
the sum of independent random variables. J Amer Statist Assoc, 1962, 57: 33--45

\bibitem{Be03} Bentkus V.  An inequality for tail
probabilities of martingales with differences bounded from one side.
J Theoret Probab, 2003, 16:  161--173

\bibitem{Be04} Bentkus V. On Hoeffding's inequality.
Ann Probab,  2004, 32:  1650--1673

\bibitem{BD14} Bentkus V, Dzindzalieta D A. Tight Gaussian bound for weighted sums of Rademacher random variables.  ArXiv:1307.3451, 2013

\bibitem{BZ06} Bentkus V,   Kalosha N,  van Zuijlen M. On domination of tail probabilities of (super)martingales: explicit bounds. Lithuanian  Math J, 2006, 46: 1--43

\bibitem{BR08} Bercu B.  In\'{e}galit\'{e}s exponentielles pour les martingales. Journ\'{e}es ALEA, 2008, 1: 1--33

\bibitem{BR06} Bercu B,    Rouault A. Sharp large deviations for the Ornstein-Uhlenbeck process. Theory Probab Appl, 2006, 46: 1--19


\bibitem{B46} Bernstein S N.  The Theory of Probabilities. Moscow, Leningrad, 1946

\bibitem{B02} Bousquet O.  A Bennett concentration inequality and its application to suprema of empirical processes. C R Acad  Sci Paris Ser I, 2002, 334:  495--500

\bibitem{CS93} Chaganty N R, Sethuraman J.   Strong large deviation and local limit theorems. Ann Probab, 1993, 21: 1671--1690

\bibitem{Cramer38} Cram\'{e}r H. Sur un nouveau th\'{e}or\`{e}me-limite de
la th\'{e}orie des probabilit\'{e}s. Actualite's Sci Indust, 1938, 736: 5--23



\bibitem{De99} de la Pe\~{n}a V H.  A general class of exponential inequalities for martingales and ratios. Ann Probab, 1999, 27: 537--564



\bibitem{D98} Dembo A,  Zeitouni O. Large deviations techniques and applications. Springer, New
York, 1998

\bibitem{D89} Deuschel J D,  Stroock D W. Large deviations. Academic Press, Boston, 1989



\bibitem{E74} Eaton M L.   A probability inequality for linear combination of bounded randon variables.  Ann Statist, 1974,  2: 609--614

\bibitem{F12} Fan X,  Grama I, Liu Q. About the constant in Talagrand's inequality for sums of bounded random variables.
ArXiv:1206.2501,  2012,  1--22

\bibitem{FGL13} Fan X, Grama I, Liu Q.  Cram\'{e}r large deviation expansions for martingales under Bernstein's condition. Stochastic Process Appl, 2013, 123: 3919--3942

\bibitem{F13} Fan X, Grama I, Liu Q.  Sharp large deviations under Bernstein's condition. C R Acad Sci Paris Ser I, 2013, 351: 845--848


\bibitem{F93} Fu J C, Li G,    Zhao L C.  On large deviation expansion of distribution of maximum likelihood estimator and its application in large sample estimation. Ann Inst Statist Math, 1993, 45: 477--498

\bibitem{IM96} It\={o} K,  MacKean H P.  Difussion Processes and Their Sample Paths.   Springer, 1996


\bibitem{GH00} Grama I, Haeusler E. Large
deviations for martingales via Cramer's method. Stochastic Process Appl, 2000, 85: 279--293

\bibitem{GHT12} Gy\"{o}rfi L,  Harrem\"{o}es P, Tusn\'{a}dy G.  Some refinements of large deviation tail probabilities.
ArXiv:1205.1005v1 [math.ST],  2012

\bibitem{Ho63} Hoeffding W.  Probability inequalities
for sums of bounded random variables. J  Amer Statist Assoc, 1963, 58: 13--30

 \bibitem{JSW03} Jing B Y, Shao Q M,  Wang  Q.  Self-normalized Cram\'{e}r-type large deviations
for independent random variables.  Ann Probab, 2003,   31: 2167--2215

\bibitem{JLZ12} Jing B Y, Liang H Y, Zhou W. Self-normalized moderate deviations for independent random variables.
Sci China Math, 2012, 55(11): 2297--2315

\bibitem{J06} Joutard C.  Sharp large deviations in nonparametric estimation. J Nonparametr Stat, 2006, 18: 293--306

\bibitem{J12} Joutard C.  Strong large deviations for arbitrary sequences of random variables. Ann Inst Stat Math, 2013, 65(1): 49--67

\bibitem{L03} Rozovky L V. A lower bound of large-deviation probabilities for the sample mean under the Cram\'{e}r condition. J  Math  Sci, 2003,  118(6)

\bibitem{L05} Rozovky L V. Large deviation probabilities for some classes of distributions statisfying the Cram\'{e}r condition. J Math Sci, 2005, 128(1)


\bibitem{N79} Nagaev S V.  Large deviations of sums
of independent random variabels. Ann Probab,  1979, 7: 745--789

\bibitem{N02} Nagaev S V.   Lower bounds for the probabilities of large deviations of sums of independent random
variables. Theory Probab  Appl, 2002,  46: 79--102; 728--735

\bibitem{Petrov75} Petrov V V.  Sums of
Independent Random Variables. Springer-Verlag, Berlin, 1975

\bibitem{Petrov95} Petrov V V.  Limit Theorems of Probability Theory.
Oxford University Press, Oxford, 1995

\bibitem{P12}  Pinelis I.  An asymptotically Gaussian bound on the Rademacher tails. Electron J Probab, 2012, 17(35):   1--22

\bibitem{P14} Pinelis I.  On the Bennett-Hoeffding inequality.  Ann Inst H Poincar\'{e} Probab Statist, 2014, 50(1):   15--27


\bibitem{R02} Rio E.  A Bennett type inequality for maxima of empirical processes.
Ann Inst H Poincar\'{e}  Probab Statist, 2002, 6: 1053--1057

\bibitem{R12q} Rio E. Sur la fonction de taux dans les in\'{e}galit\'{e}s de Talagrand pour les processus empiriques.
C R Acad Sci Paris Ser I, 2012, 350: 303--305

\bibitem{R12} Rio E. In\'{e}galit\'{e}s exponentielles et in\'{e}galit\'{e}s de concentration. Institut Math\'{e}matique de Bordeaux, 2012, 1: 1--22.

\bibitem{SS91} Saulis L, Statulevicius V A. Limit theorems for large deviations. Springer, 1991


\bibitem{S10} Shevtsova I G.  An improvement of convergence rate estimates
in the Lyapunov theorem. Doklady  Math, 2010, 82: 862--864

\bibitem{S91} Sakhanenko A I.  Berry-Esseen type bounds for large deviation probabilities. Siberian Math J, 1991, 32: 647--656



\bibitem{S99} Shao Q M.  A Cram\'{e}r type large deviation result for Student's t-statistic.  J Theoret Probab,   1999, 12(2): 385--398

\bibitem{Ta95} Talagrand M.  The missing factor in
Hoeffding's inequalities. Ann Inst  H  Poincar\'{e} Probab Statist, 1995, 31: 689--702


\bibitem{VL95} van de Geer S.
Exponential inequalities for martingales, with application to maximum likelihood estimation for counting process.
Ann Statist, 1995, 23: 1779--1801

\bibitem{VL11} van de Geer S, Lederer J.
The Bernstein-Orlicz norm and deviation inequalities. Probab Theory Relat Fields,  2013, 157:  225--250



\end{thebibliography}
\end{document}